\documentclass[leqno]{article}
\usepackage{amssymb,latexsym,amscd}
\newcommand{\se}[1]{{\section{#1}} {\setcounter{equation}{0}}}
\newtheorem{theorem}{Theorem}[section]
\newtheorem{lm}{Lemma}[section]
\newtheorem{prop}{Proposition}[section]
\newtheorem{de}{Definition}[section]
\newtheorem{co}{Corollary}[section]

\def\k{{K\"{a}hler }}

\def\cy{{Calabi-Yau }}
\def\l{{Lagrangian }}

\input epsf
\begin{document}
\hbadness=10000
\title{{\bf Lagrangian torus fibration of quintic Calabi-Yau hypersurfaces III:}\\
{\large {\bf Symplectic topological SYZ mirror construction for general quintics}}}
\author{Wei-Dong Ruan\\
Department of Mathematics\\
University of Illinois at Chicago\\
Chicago, IL 60607\\}
\date{}
\footnotetext{Partially supported by NSF Grant DMS-9703870 and DMS-0104150.}
\maketitle
\begin{abstract}
In this article we construct Lagrangian torus fibrations for general quintic \cy hypersurfaces near the large complex limit and their mirror manifolds using gradient flow method. Then we prove the Strominger-Yau-Zaslow mirror conjecture for this class of \cy manifolds in symplectic category.
\end{abstract}
\se{Introduction}
This paper is a sequel of \cite{lag1, lag2}. The motivation of studying Lagrangian torus fibrations of \cy manifolds comes from Strominger-Yau-Zaslow's proposed approach toward mirror symmetry (\cite{SYZ}). According to their proposal, on each Calabi-Yau manifold, there should be a special Lagrangian torus fibration. This conjectural special Lagrangian torus fibration structure of Calabi-Yau manifolds is further used to give a possible geometric explanation of mirror symmetry conjecture. Despite its great potential in solving the mirror symmetry conjecture, there are very few known examples of special Lagrangian submanifolds or special Lagrangian fibrations for dimension $n\geq 3$. Given our lack of knowledge for special Lagrangian, one may consider relaxing the requirement and consider Lagrangian fibration, which is largely unexplored and interesting in its own right. For many applications to mirror symmetry, especially those concerning (symplectic) topological structure of fibrations, Lagrangian fibrations will provide quite sufficient information. In this paper, as in the previous two papers (\cite{lag1}, \cite{lag2}), we will mainly concern Lagrangian torus fibrations of Calabi-Yau hypersurfaces in toric variety, namely, the symplectic topological aspect of SYZ mirror construction.\\

In this paper, we will construct Lagrangian torus fibrations of generic quintic Calabi-Yau hypersurfaces and their mirror manifolds in complete generality. With the detailed understanding of Lagrangian torus fibrations, we will be able to prove the symplectic topological SYZ mirror conjecture for Calabi-Yau quintic hypersurfaces in $\mathbb{CP}^4$. More precisely,\\

\begin{theorem}
For generic quintic Calabi-Yau hypersurface $X$ near the large complex limit, and its mirror Calabi-Yau manifold $Y$ near the large radius limit, there exist corresponding Lagrangian torus fibrations\\
\[
\begin{array}{ccccccc}
X_{s(b)}&\hookrightarrow& X& \ \ &Y_b&\hookrightarrow& Y\\
&&\downarrow& \ &&&\downarrow\\
&& \partial \Delta& \ &&& \partial \Delta_w\\
\end{array}
\]\\
with singular locus $\Gamma \subset \partial \Delta$ and $\Gamma' \subset \partial \Delta_w$, where $s:\partial \Delta_w \rightarrow \partial \Delta$ is a natural homeomorphism, $s(\Gamma')=\Gamma$. For $b\in \partial \Delta_w \backslash \Gamma'$, the corresponding fibers $X_{s(b)}$ and $Y_b$ are naturally dual to each other.
\end{theorem}

This theorem will be proved in more precise form as in theorem \ref{eb}.\\

{\bf Remark:} It is important to point out that our purpose is not merely to construct {\bf a} Lagrangian fibration for Calabi-Yau manifold. We want to construct {\bf the} Lagrangian fibration corresponding to the K\"{a}hler and the complex moduli of the Calabi-Yau manifold. Namely, for Calabi-Yau manifolds with different K\"{a}hler and complex structures, the structures of the corresponding Lagrangian fibration (singular set, singular locus, singular fibers, etc.) should be different and reflect the corresponding K\"{a}hler and complex structures.\\

The gradient flow approach we developed in \cite{lag1, lag2} is ideal for this purpose. In addition to its clear advantage of being able to naturally produce {\bf Lagrangian} fibrations, for general Calabi-Yau hypersurfaces, the gradient flow will naturally produce a rather ``canonical'' family of Lagrangian torus fibrations continuously depending on the moduli of complex and K\"{a}hler structures of the Calabi-Yau hypersurface. Some part of the fibration structure (singular set, singular locus, singular fibers, etc.) will depend on the complex structure, and some other part will depend on the K\"{a}hler structure. It is important to understand the dependence precisely. For general Calabi-Yau hypersurfaces, the dependence on complex and K\"{a}hler structure is somewhat mixed. For quintics, the complex moduli has 101 dimensions and the K\"{a}hler moduli has only 1 dimension, so the fibration structures mainly depend on the complex moduli. On the other side, for the mirrors of quintics, the K\"{a}hler moduli has 101 dimensions and the complex moduli has only 1 dimension, so the fibration structures mainly depend on the K\"{a}hler moduli. Therefore, it is worthwhile to first discuss the case of quintics and their mirrors in detail, where the dependence can be seen more clearly. This will help us better understand the dependence of the fibration structures on complex and K\"{a}hler moduli for more general cases. Take singular locus as an example. According to the discussion in section 2, the structure of singular locus for the quintics is largely a result of different ``string diagram'' structures of the singular set curves near the large complex limit. According to the discussion in section 4, the structure of singular locus for the mirrors of quintics is largely a result of the combinatorial structures of different crepant resolutions of singularities. For more general cases that we will discuss in \cite{tor,ci}, both phenomena will come into effect together in forming the singular locus.\\

As a byproduct, the detailed understanding of the Lagrangian torus fibration structures also provides better understanding of the mirror symmetry inspired partial compactification of the complex moduli of Calabi-Yau manifolds. As is well known, for the sake of mirror symmetry, the complex moduli of Calabi-Yau hypersurfaces in toric varieties near the large complex limit should be partially compactified according to the secondary fan (see \cite{AGM}, \cite{AGM2}, etc.). Authors of \cite{AGM,AGM2} also proposed a related chamber decomposition of the complex moduli of Calabi-Yau hypersurfaces near the large complex limit so that the chambers under the monomial-divisor mirror map correspond to the \k cones of the different birationally equivalent models of the mirror, which make up the top dimensional cones of the secondary fan. To our knowledge, so far there is no satisfactory intrinsic geometric explanation of this chamber decomposition nor convincing direct reasoning (without going to the \k moduli of the mirror) why secondary fan compactification should be the suitable partial compactification for mirror symmetry purpose. Our construction of Lagrangian torus fibrations as application gives an intrinsic geometric explanation of this chamber decomposition, consequently provides direct reasoning for the secondary fan compactification. More precisely, the Calabi-Yau manifolds in the same chamber near the large complex limit are exactly those whose Lagrangian torus fibrations have the graph singular locus with the same combinatorial type. This interpretation (discussed in section 3.1) can potentially be generalized to determine the chamber decomposition and partial compactification of the complex moduli near the large complex limit for more general Calabi-Yau manifolds in non-toric situations.\\

Next we give a brief review of our work in the first two papers \cite{lag1,lag2}.\\

Our work on Lagrangian torus fibrations of Calabi-Yau manifolds starts with \cite{lag1}. In \cite{lag1}, we described a very simple and natural construction of Lagrangian torus fibrations via gradient flow based on a natural Lagrangian torus fibration of the large complex limit. This method in principle will produce Lagrangian torus fibrations for general Calabi-Yau hypersurfaces in toric varieties. For simplicity, we described the case of Fermat type quintic Calabi-Yau threefold family $\{X_\psi\}$ in $\mathbb{CP}^4$ defined by:

\[
p_{\psi}=\sum_1^5 z_k^5 - 5\psi \prod_{k=1}^5 z_k=0
\]

near the large complex limit $X_{\infty}$

\[
p_{\infty}=\prod_{k=1}^5 z_k=0
\]

in great detail. Most of the essential features of the gradient flow for the general cases already show up there. One crucial discovery in \cite{lag1} is that although the Lagrangian torus fibration of $X_{\infty}$ has many lower-dimensional torus fibers, the gradient flow will automatically yield Lagrangian torus fibration for smooth $X_\psi$ with 3-dimensional fibers everywhere and singular fibers clearly located. We also discussed the ``expected" structure of the special Lagrangian torus fibration, particularly we computed the monodromy transformations of the ``expected" special Lagrangian fibration and discussed the ``expected" singular fiber structures implied by the monodromy information. Then we compared our Lagrangian fibration constructed via gradient flow with the ``expected" special Lagrangian fibration and noted the differences. Finally, we discussed its relevance to mirror construction for Calabi-Yau hypersurfaces in toric varieties.\\

The so-called ``expected" special Lagrangian fibrations first used in \cite{lag1} refers to the generally expected structure of the special Lagrangian torus fibrations of Calabi-Yau manifold at the time of \cite{lag1}, based on our knowledge of the elliptic fibrations of K3 surfaces corresponding to SYZ in dimension 2. More precisely, the special Lagrangian fibration was expected to be $C^\infty$, which necessarily has codimension 2 singular locus. The codimension 2 singular locus condition is crucial for SYZ conjecture as originally proposed to be valid for special Lagrangian fibrations. The Lagrangian fibrations we constructed in \cite{lag1}, with  codimension 1 singular locus and different singular fibers, do not have exactly the ``expected" structure of a $C^\infty$ special Lagrangian fibration. Recent examples of D. Joyce \cite{Joyce} indicate that the structure of the actual special Lagrangian fibrations are probably not as expected, instead, they probably resemble more closely the structure of the natural Lagrangian fibrations we constructed in \cite{lag1} with codimension 1 singular locus. (This suggests that the singular loci of the special Lagrangian fibrations for a Calabi-Yau manifold and its mirror are probably different, based on the structure of codimension 1 singular loci of our Lagrangian fibrations for quintic Calabi-Yau and its mirror, which will be discussed in detail in this paper. Consequently, the original SYZ conjecture for special Lagrangian fibrations need to be modified.) This recent development makes the Lagrangian torus fibrations with codimension 1 singular locus potentially more important as a starting point to deform to the actual special Lagrangian fibration. Nevertheless, it is still interesting and important to construct and study Lagrangian torus fibrations with codimension 2 singular locus, because firstly we can show the symplectic SYZ mirror conjecture (as precisely formulated in section 6) to be valid for this kind of fibrations, secondly they have many good properties and are very useful for topological computations, and also hold independent interest from symplectic geometric point of view. In this paper, the term ``expected" Lagrangian torus fibrations will more narrowly refer to the topological structure of Lagrangian torus fibrations with codimension 2 topological singular locus.\\

From the above discussion, special Lagrangian fibrations in SYZ construction are likely to be non-smooth. Coincidentally, our gradient flow approach naturally produces piecewise smooth (Lipschitz) Lagrangian fibrations. However, in our opinion, $C^\infty$ Lagrangian fibrations should still play an important role for symplectic topological aspect of SYZ construction. This is indicated by a somewhat surprising fact discussed in \cite{lag2}, that a general piecewise smooth (Lipschitz) Lagrangian fibration can not be smoothed to a $C^\infty$ Lagrangian fibration by small perturbation. More precisely, the singular locus of a general piecewise smooth (Lipschitz) Lagrangian fibration is usually of codimension 1 while the singular locus of the corresponding $C^\infty$ Lagrangian fibration necessarily has codimension 2. For this reason, it is desirable to try to find out how smooth one can make the naturally constructed piecewise smooth (Lipschitz) Lagrangian fibrations to be. When the smoothing can be achieved, we get a nice fibration with well behaved singular locus and singular fibers. When the smoothing can only be partially achieved, the obstruction will give us better understanding of the symplectic topology of the Calabi-Yau manifold and shed a light on the reason why the special Lagrangian fibration can not be smooth in general. In \cite{lag2}, we constructed many local examples of piecewise smooth (Lipschitz) Lagrangian fibrations that are not $C^\infty$. We also discussed methods to squeeze the codimension 1 singular locus of our piecewise smooth (Lipschitz) Lagrangian fibration into codimension 2 by symplectic geometry techniques. This can be viewed as the first step (topological modification) toward the smoothing of our Lagrangian torus fibration. Partial smoothing (analytical modification) of such Lagrangian fibration with codimension 2 singular locus will be discussed in \cite{smooth}.\\

For the understanding of the topological aspect of SYZ, a construction of purely topological (non-Lagrangian) torus fibration can also be of interest. Our method clearly can be easily modified to construct non-Lagrangian torus fibrations with much less technical difficulty and better smoothness. We choose to construct more difficult Lagrangian fibrations because first of all our gradient flow naturally produces Lagrangian fibrations. Secondly, Lagrangian fibrations could impose strong constraint on the topological types of singular fibers (see section 2 of \cite{lag2}). The very difficulty involved in smoothing the Lagrangian fibrations can be taken as an indication that the actual special Lagrangian fibrations may not be smooth in general ( also indicated by recent examples of D. Joyce \cite{Joyce}). Non-Lagrangian fibrations, which can easily be smoothed, would not give us such insight. In light of all these, Lagrangian fibrations seem to be a good compromise between the very rigid special Lagrangian fibrations and general non-Lagrangian torus fibrations, which lack control of singular fibers.\\

In \cite{lag2} we provided the technical details for the gradient flow construction. The flow we use to produce Lagrangian fibration is the so-called normalized gradient flow. Gradient flows of smooth functions with non-degenerate critical points have been used intensively in the mathematical literature, for example in Morse theory. But the gradient flow in our situation is very unconventional. The critical points of our function are usually highly degenerate and often non-isolated. Worst of all our function is not even smooth (it has infinities along some subvarieties). In \cite{lag2}, we computed local models for the singular vector fields, then proved a structural stability for the gradient vector fields we use, which ensures that the gradient flow would behave as we expect. Other technical aspects in \cite{lag2} include the deformation of symplectic manifold and sympletic submanifold structures, the construction of toroidal \k metrics, and the symplectic deformation of the Lagrangian torus fibration with codimension 1 singular locus into a Lagrangian torus fibration with codimension 2 singular locus.\\

In the following we briefly discuss the contents of this paper.\\

As we know mirror symmetry conjecture first of all amounts to identifying the complex moduli of a Calabi-Yau manifold with the complexified K\"{a}hler moduli of its mirror manifold, namely, specifying the mirror map. Then it further concludes that the quantum geometry of a Calabi-Yau manifold is equivalent to the quantum geometry of its mirror partner. Despite the important role of Fermat type quintic family in the history of mirror symmetry, and the nice features of Lagrangian fibrations and singular fibers of this family as discussed in \cite{lag1} and \cite{lag2}, the Fermat type quintic Calabi-Yau family is a highly symmetric and very special type of Calabi-Yau manifolds. In fact, this family is located at the boundary of the moduli space of quintics. Their mirror manifolds are highly degenerate Calabi-Yau orbifolds that are located at the boundary of the mirror K\"{a}hler cone. To understand mirror symmetry for quintic Calabi-Yau threefolds, it is more important to understand Lagrangian torus fibrations of {\bf generic} quintic Calabi-Yau hypersurfaces in $\mathbb{CP}^4$ and their mirrors that we will discuss in this paper.\\

In Section 2 we construct Lagrangian torus fibrations for generic quintics, first with codimension 1 singular locus, then with codimension 2 singular locus. The construction of Lagrangian torus fibrations with codimension 1 singular locus is automatic using the gradient flow method, just as for the Fermat quintics, and works the same way even for more general situations. The main difficulties come from understanding the structure of the resulting codimension 1 singular locus and squeezing the codimension 1 singular locus into codimension 2. For Fermat type quintics, the codimension 1 singular locus can be easily seen as a fattening of a graph, so we can always modify the Lagrangian fibrations with codimension 1 singular locus to the expected topological type. It seems hopeless to do the same for general quintics due to the generally badly behaved singular locus. It turns out that for general quintics {\bf near} the large complex limit (in suitable sense), the codimension 1 singular locus is much better behaved, and we can do the same thing as in the Fermat type Calabi-Yau case to get the expected Lagrangian torus fibrations.\\

More precisely, let $F: X\rightarrow S^3$ be the Lagrangian torus fibration of a quintic constructed by the gradient flow. The singular set $C\subset X$ is a complex curve (more precisely, $C={\rm Sing}(X_\infty)\cap X$). The singular locus $\tilde{\Gamma}=F(C)$ is usually a 2-dimensional object (called amoeba in \cite{G}), which in general can be rather chaotic. Miraculously, when $X$ is near the large complex limit (in a suitable sense defined in section 2), according to results in our paper \cite{N}, the singular locus $\tilde{\Gamma}$ is actually a fattening of some 1-dimensional graph $\Gamma$. With this fact and corresponding squeezing result in \cite{N}, general methods developed in \cite{lag1} and \cite{lag2} will enable us to construct Lagrangian torus fibrations with codimension 2 graph singular locus for generic quintics near the large complex limit.\\

In section 3 we discuss the complex moduli for the quintics and the \k moduli for their mirror manifolds. In subsection 3.1, we provide a direct geometric reconstruction of the secondary fan compactification of the complex moduli for the quintics based on our construction of Lagrangian torus fibrations for Calabi-Yau quintics from section 2. The discussion of the SYZ mirror construction needs an identification of the complex moduli of a Calabi-Yau near the large complex limit with the complexified K\"{a}hler moduli of the mirror Calabi-Yau to start with. According to the work of Aspinwall-Greene-Morrison in \cite{AGM}, for Calabi-Yau hypersurfaces in toric variety, it is more natural to consider the monomial-divisor mirror map instead of the actural mirror map, which is a higher order perturbation of the monomial-divisor mirror map. Our symplectic topological SYZ construction will be based on the monomial-divisor mirror map. Results in \cite{AGM} do not directly apply to our case. By using a slicing theorem we remove some restrictions in \cite{AGM} and construct the monomial-divisor mirror map for quintic Calabi-Yau hypersurfaces in the form we need based on \cite{AGM}.\\

Lagrangian torus fibrations of the mirrors of quintics are constructed in section 4. As proposed by Greene, Plesser and others (see \cite{GP, Roan}), the mirror manifolds of quintics can be understood as crepant resolutions of orbifold quotient of Fermat type quintics. The Lagrangian torus fibration structures of the mirror of quintics mostly depend on the \k moduli, which can be reduced to the classical knowledge of resolution of singularities. Indeed, most part of the singular locus of our natural Lagrangian torus fibration for the mirror of quintic with codimension 1 singular locus is actually of codimension 2, which has a nice interpretation in terms of resolution of singularities. The remaining codimension 1 part of singular locus is a disjoint union of finitely many fattened ``Y". The discussion here on squeezing to codimension 2 singular locus is a bit tricky, because locally the \k metric may not be Fubini-Study. Besides the discussion of construction of Lagrangian torus fibrations for the mirrors of quintics with codimension 1 and codimension 2 singular locus (in subsections 5.1 and 5.3) we also discuss the construction of non-Lagrangian torus fibration with codimension 2 singular locus (in subsection 5.2), which is much easier.\\

The actual symplectic topological SYZ construction is worked out in section 5. We first identify the bases and the singular locus of the two fibrations, then we establish the duality of the regular fibers. In the process, we also find a very simple way to compute the monodromy of the Lagrangian torus fibrations. In section 6, certain classes of generic singular fibers of Lagrangian torus fibrations are discussed and a proposal of SYZ mirror duality for such generic singular fibers is presented.\\

The gradient flow method and the construction of Lagrangian torus fibrations, slicing theorem, monomial-divisor map, construction of singular locus, identification of base spaces, duality relation of the fibers and monodromy computation, etc., can all be generalized to the general situation of Calabi-Yau hypersurfaces in a toric variety corresponding to a reflexive polyhedron based on the methods and constructions we developed in \cite{lag1}, \cite{lag2}, \cite{N} and in this paper. These generalizations and the symplectic topological SYZ mirror conjecture for general generic Calabi-Yau hypersurfaces in a toric variety corresponding to a reflexive polyhedron is proved in complete generality in \cite{tor}. In \cite{ci}, our construction is further generalized to Calabi-Yau complete intersections in toric varieties.\\

There is an earlier work by Zharkov (\cite{Z}) on the construction of certain non-Lagrangian torus fibrations of Calabi-Yau hypersurfaces in toric variety. The work of Leung and Vafa (\cite{lv}) from physics point of view, touched upon several important ideas related to our work. There is also the work of M. Gross \cite{GW}, which appeared on the internet at around the same time as this paper, where certain non-Lagrangian torus fibrations for quintic Calabi-Yau are constructed based on the information of the torus fibration structure for the mirror of quintic and with the help of C.T.C. Wall's existence theorem. One major difference between our work and these other results is that our torus fibrations are naturally Lagrangian fibrations.\\

{\bf Note on notation:} Unless otherwise specified, Lagrangian fibration in this paper only refers to piecewise smooth (Lipschitz) Lagrangian fibration. The notion of convex in this paper is sometimes non-standard. In certain cases, it probably is called concave in conventional term. For precise definition, please refer to the respective sections. Another thing is that we usually use $\Delta$ to denote the set of integral points in a Newton polyhedron, but sometimes we use $\Delta$ to denote the corresponding real polyhedron.\\

\se{Lagrangian torus fibrations of general quintic Calabi-Yau threefolds}
In \cite{lag1}, \cite{lag2}, we constructed Lagrangian torus fibrations for Fermat type quintic Calabi-Yau family $\{X_\psi\}$ in $\mathbb{CP}^4$ by a flow along vector fields. It is easy to observe that the same construction will also produce Lagrangian torus fibrations for general quintics. Let\\
\[
z^m=\prod_{k=1}^5z_k^{m_k},\ |m|=\sum_{k=1}^5 m_k,\ {\rm for}\ m=(m_1,m_2,m_3,m_4,m_5)\in \mathbb{Z}^5_{\geq 0}.
\]
Then a general quintic homogeneous polynomial can be written as\\
\[
p(z) = \sum_{|m|=5}a_mz^m.
\]\\
Let $m_0=(1,1,1,1,1)$, and denote $a_{m_0}=\psi$. Consider the quintic Calabi-Yau family $\{X_\psi\}$ in $\mathbb{CP}^4$ defined by\\
\[
p_\psi(z) = p_a(z) + \psi \prod_{k=1}^5 z_k = \sum_{m\not=m_0,|m|=5}a_mz^m + \psi \prod_{k=1}^5 z_k=0.
\]\\
When $\psi$ approaches $\infty$, the family approaches its ``large complex limit'' $X_{\infty}$ defined by\\
\[
p_{\infty}=\prod_{k=1}^5 z_k=0.
\]\\
Consider the meromorphic function\\
\[
s= \frac{p_\infty(z)}{p_a(z)}
\]\\
defined on $\mathbb{CP}^4$. Let $\omega$ denote the K\"{a}hler form of a K\"{a}hler metric $g$ on $\mathbb{CP}^4$, and $\nabla f$ denote the gradient vector field of the real function $f=Re(s)$ with respect to the K\"{a}hler metric $g$. As in \cite{lag1}, we will similarly use the flow of the normalized gradient vector field $V=\frac{\nabla f}{|\nabla f|^2}$ to construct Lagrangian torus fibration for $X_\psi$ based on Lagrangian torus fibration of $X_\infty$.\\

From our experience on Fermat type quintic Calabi-Yau hypersurfaces, we know that the Lagrangian fibration we get by this flow usually has codimension 1 singular locus. Due to the highly symmetric nature of the Fermat type quintics, in \cite{lag1} it was relatively easy to figure out explicitly the ``expected" Lagrangian fibration structure with codimension 2 singular locus and the corresponding singular fibers. In that case, the singular locus of the Lagrangian fibration with codimension 1 singular locus is a fattened version of the codimension 2 singular locus of the ``expected" Lagrangian fibration. Therefore, even without perturbing to the Lagrangian fibration with codimension 2 singular locus, we could already compute the monodromy. A rather magical explicit symplectic deformation construction (discussed in section 9 of \cite{lag2}) was used to deform the Lagrangian torus fibration with codimension 1 singular locus to one with codimension 2 singular locus.\\

The major difficulty to generalize this program to general quintic Calabi-Yau threefolds is that for general quintics, the singular locus of the Lagrangian fibration for $X_\psi$ constructed from deforming the standard Lagrangian torus fibration of $X_{\infty}$ via the flow of $V$ can be fairly arbitrary and does not necessarily resemble the fattening of any ``expected singular locus''. Worst of all, in the case of general quintic Calabi-Yau threefolds, there is no obvious guess what the ``expected" codimension 2 singular locus should be and how it will vary when deforming in the complex moduli of the quintic Calabi-Yau. One clearly expects the ``expected" singular locus to be some graph in $\partial \Delta \cong S^3$. But it seems to take some miracle (at least to me when I first dreamed about it) for a general singular set $C$ (which is an algebraic curve) to project to the singular locus $\tilde{\Gamma} = F(C)$ that resembles a fattening of a graph.\\

Interestingly, miracle happens here! It largely relies on a better understanding of what it means to be {\bf near the large complex limit}. The large complex limit was a rather elusive concept from physics. From early works of \cite{Roan, GP, Can}, it is apparent that $X_\infty$ can be viewed in some sense as a representative of the large complex limit for quintic Calabi-Yau. (We will often refer to $X_\infty$ in this paper as the large complex limit in this narrow sense.) Later, authors of \cite{AGM2} and others proposed that there are many large complex limits corresponding to a chamber decomposition of the complex moduli of quintic Calabi-Yau near $X_\infty$. (Each chamber under the monomial-divisor mirror map corresponds to the complexified \k cone of one of the birationally equivalent model of the mirror of quintic.) Each chamber characterize a particular way quintics get near $X_\infty$. This can be interpreted equivalently from two different perspectives as either different ways quintics get near the same large complex limit or quintics approaching different large complex limits. (In this paper, we will more often adopt the first interpretation.) For our purpose, when quintics get deep into each such chamber and approach $X_\infty$, the singular locus of our fibration will resemble fattening of a graph. For different chambers the corresponding graphs will be different. From this point of view, our discussion of singular locus in terms of amoeba actually will give a very natural new way to reconstruct this chamber decomposition without going to the mirror, which will be discussed in more detail in section 3 when we talk about the complex moduli space of quintic.\\

Before getting into the detail, let us recall the general result on the construction of Lagrangian torus fibration (theorem 8.1) in \cite{lag2}. In the following, we will rephrase this general theorem according to our special situation here. Assume that

\[
X_{\rm inv} = X_\psi \cap X_\infty,\ \ {\rm in\ particular}\ \ C = X_\psi \cap {\rm Sing}(X_\infty).
\]

For any subset $I \subset \{1,2,3,4,5\}$, let

\[
D_I = \{z\in \mathbb{CP}^4| z_i=0,z_j\not= 0,\ {\rm for}\ i\in I,j\in\{1,2,3,4,5\}\backslash I\}.
\]

In this section, we will {\bf always} use $I$ to denote a subset in $\{1,2,3,4,5\}$. Let $|I|$ denote the cardinality of $I$, then

\[
X_\infty = \bigcup_{\tiny{\begin{array}{c}I \subset \{1,2,3,4,5\}\\ 0<|I|<5\end{array}}} D_I.
\]

Let $\Delta$ denote the standard 4-simplex, whose vertices are identified with $\{1,2,3,4,5\}$. We will use $\Delta_I$ to denote the subface of $\Delta$, whose vertices are not in $I$. Then we have

\[
\partial\Delta = \bigcup_{\tiny{\begin{array}{c}I \subset \{1,2,3,4,5\}\\ 0<|I|<5\end{array}}} \Delta_I.
\]

\begin{de}
Assume that $(M_1,\omega_1)$ is a smooth symplectic manifold and $(M_2,\omega_2)$ is a symplectic variety. Then a piecewise smooth map $H: M_1 \rightarrow M_2$ is called a symplectic morphism if $H^*\omega_2 = \omega_1$. If $(M_2,\omega_2)$ is also a smooth symplectic manifold and $H$ is a diffeomorphism, then the symplectic morphism $H$ is also called a symplectomorphism.
\end{de}

{\bf Remark:} In this paper, we will only deal with the case of normal crossing symplectic varieties. More specifically, in our case, $M_2 = X_\infty$ and $\omega_2$ is taken to be the restriction to $X_\infty$ of a symplectic form on $\mathbb{CP}^4$. Therefore, we will not venture into the concept of general symplectic varieties and symplectic forms on them.\\

\begin{de}
When $X$ and $B$ are smooth, a fibration $F: X \rightarrow B$ is called topologically smooth if $F$ is locally a topological product. When $X$ and $B$ are stratified into unions of smooth strata

\[
X = \bigcup_i X_i,\ \ B = \bigcup_i B_i,
\]

such that $F(X_i) = B_i$, $F: X \rightarrow B$ is called topologically smooth (with respect to the stratifications) if $F$ is continuous and $F|_{X_i}: X_i \rightarrow B_i$ are topologically smooth for all $i$.\\
\end{de}

Notice that the topological smoothness of $F$ depends on the stratifications of $X$ and $B$. Specifying to our case of torus fibration, a torus fibration $F_\infty: X_\infty \rightarrow \partial \Delta$ is topologically smooth if for each $I$, $\Delta_I = F_\infty(D_I)$ and $F_\infty|_{D_I}: D_I \rightarrow \Delta_I$ is a topologically smooth fibration, whose fibers are necessarily $(4-|I|)$-dimensional torus. Using our notation here, the theorem 8.1 of \cite{lag2} can be rephased in our situation as follows.\\

\begin{theorem}
\label{bf}
Start with a topologically smooth Lagrangian torus fibration $F_\infty: X_\infty \rightarrow \partial \Delta$, we can construct a symplectic morphism $H_\psi: X_\psi \rightarrow X_\infty$ such that $F_\psi = F_\infty\circ H_\psi: X_\psi \rightarrow \partial \Delta$ is a Lagrangian torus fibration with singular set $C = X_\psi \cap {\rm Sing}(X_\infty)$ and singular locus $\Gamma = F_\infty(C)$. For $b\not\in \Gamma$, $F_\psi^{-1}(b)$ is a real $3$-torus. For $b \in \Gamma$, $F_\psi^{-1}(b)$ is singular. For $b \in \Gamma \cap \Delta_I$, $F_\psi^{-1}(b)\cap C = F_\infty^{-1}(b)\cap C$ and $H_\psi: F_\psi^{-1}(b)\backslash C \rightarrow F_\infty^{-1}(b)\backslash C$ is a topologically smooth $(|I|-1)$-torus fibration.
\end{theorem}
\begin{flushright} $\Box$ \end{flushright}
{\bf Remark:} Since our map $F_\psi$ is not $C^\infty$, by saying a point is in the singular set of $F_\psi$, we mean that in a neighborhood of that point the fibration $F_\psi$ is not a topological product.\\

{\bf Remark on notation:} We generally use $\Gamma$ to denote the singular locus. Sometimes when we discuss the Lagrangian torus fibrations with codimension 1 singular locus and the related Lagrangian torus fibrations with codimension 2 singular locus together, we usually use $\tilde{\Gamma}$ to denote the codimension 1 singular locus and $\Gamma$ to denote the codimension 2 singular locus.\\

From theorem \ref{bf}, the construction of Lagrangian fibration for $X_\psi$ can be reduced to the construction of the topologically smooth Lagrangian torus fibration $F_\infty: X_\infty \rightarrow \partial \Delta$. The singular locus $\Gamma = F_\infty(C)$ is determined by $F_\infty$ and $C$. The singular fibers are determined by the ways fibers of $F_\infty$ intersect with $C$. For a fixed quintic family $\{X_\psi\}$, $C$ is fixed. Therefore from now on our discussion will mainly be focused on the construction of $F_\infty$ in various situations.\\

The most natural choice for $F_\infty$ is the restriction to $X_\infty$ of the moment map $F_{\rm FS}: \mathbb{CP}^4 \rightarrow \Delta$ with respect to the Fubini-Study \k form. The fibers of $F_\infty = F_{\rm FS}|_{X_\infty}: X_\infty \rightarrow \partial \Delta$ are exactly the orbits of the real torus action. For this reason, we may use any toric \k forms on $\mathbb{CP}^4$. They will give the same torus fibration for $X_\infty$ except for a difference of reparametrization of the base $\partial \Delta$. $F_\infty$ is clearly a topologically smooth Lagrangian torus fibration. By theorem \ref{bf}, we can construct a Lagrangian torus fibration $F_\psi: X_\psi \rightarrow \partial \Delta$ with the singular locus $\tilde{\Gamma} = F_\infty(C)$.\\

For $I \subset \{1,2,3,4,5\}$ with $|I|=2$, let

\[
C_I =\{[z]\in \mathbb{CP}^4| p_a(z) =0, z_l=0\ \ {\rm for}\ l\in I\}.
\]

$C_I$ is a genus 6 curve. We have

\[
C = {\rm Sing}(X_{\infty})\cap X_{\psi} = \bigcup_{\tiny{\begin{array}{c}I \subset \{1,2,3,4,5\}\\ |I|=2\end{array}}} C_I.
\]

The singular locus

\[
\tilde{\Gamma} = F_\infty(C) = \bigcup_{\tiny{\begin{array}{c}I \subset \{1,2,3,4,5\}\\ |I|=2\end{array}}} \tilde{\Gamma}_I = \bigcup_{\tiny{\begin{array}{c}I \subset \{1,2,3,4,5\}\\ |I|=2\end{array}}} F_\infty(C_I).
\]

At this point, we already constructed a Lagrangian fibration $F_\psi: X_\psi \rightarrow \partial \Delta$ with the codimension 1 singular locus $\tilde{\Gamma} = F_\infty(C)$. To connect with the SYZ picture, especially to construct Lagrangian torus fibration with graph singular locus, more detailed understanding of the structure of the singular locus $\tilde{\Gamma} = F_\infty(C)$ is needed.\\

What we need is to realize the singular locus $\tilde{\Gamma} = F_\infty(C)$ as a fattening of some graph $\Gamma$ and be able to explicitly construct a symplectic deformation that deforms the \l fibration $F_\infty$ to $\hat{F}$ such that $\hat{F}(C) = \Gamma$. According to the section 9 of \cite{lag2}, this is equivalent to symplectically deforming $C$ to a symplectic curve $\hat{C}$ that satisfies $F_\infty(\hat{C})=\Gamma$. Since $C$ is reducible, and each irreducible component $C_I$ is in $\overline{D_I} \cong \mathbb{CP}^2$, similar to the argument in the section 9 of \cite{lag2}, the problem can be isolated to each $\mathbb{CP}^2$ and be reduced to the following problem. We will temporarily use $\Delta$ to denote a 2-simplex in the following.\\

{\bf Problem:} Let $F:\mathbb{CP}^2 \rightarrow \Delta$ be the standard moment map with respect to certain toric metric on $\mathbb{CP}^2$. We need to find quintic curves $C$ in $\mathbb{CP}^2$, such that $\tilde{\Gamma} = F(C)$ is a fattening of some graph $\Gamma$. We also want to explicitly construct a symplectic deformation of $\mathbb{CP}^2$, which deforms $C$ to symplectic curve $\hat{C}$ with $F(\hat{C})=\Gamma$.\\

Clearly, one can not expect every quintic curve to have such nice properties. As we mentioned earlier, it turns out that when the Calabi-Yau quintic is generic and close to the large complex limit in a certain sense, the corresponding quintic curves $C_{ijk}$ have the properties described above. This kind of curves and more general situations have been discussed intensively in our paper \cite{N}. To describe the result from that paper, we first introduce some notations.\\

Consider $\mathbb{CP}^2$ with homogeneous coordinate $[z_0,z_1,z_2]$ and inhomogeneous coordinate $(x_1,x_2)=(z_1/z_0,z_2/z_0)$. Let

\[
M = \{ x^m = x_1^{m_1}x_2^{m_2} | m=(m_1,m_2)\in \mathbb{Z}^2\} \cong \mathbb{Z}^2.
\]

Consider a general quintic polynomial

\[
p(x)= \sum_{|m|\leq 5}a_m x^m,
\]

where $|m|=m_1+m_2$. Then the Newton polygon of quintic polynomials is

\[
\Delta_5 = \{m\in M | m_1,m_2 \geq 0, |m|\leq 5\}.
\]

Let $w=(w_m)_{m\in \Delta_5}$ be a function on $\Delta_5$ (regarding $\Delta_5$ as an integral triangle in the lattice $M$).\\
\begin{de}
\label{ba}
$w=(w_m)_{m\in \Delta_5}$ is called convex on $\Delta_5$ if for any $m' \in \Delta_5$, there exists an affine function $n$ such that $n(m')= w_{m'}$ and $n(m)\leq w_m$ for $m\in \Delta_5\backslash \{m'\}$.\\
\end{de}
We will always assume that $w$ is convex. With $w$ we can define the moment map\\
\[
F_{t^w}(x) = \sum_{m\in\Delta_5} \frac{|x^m|_{t^w}^2}{|x|_{t^w}^2}m,
\]\\
which maps $\mathbb{CP}^2$ to $\Delta_5$, where $t>0$ is a parameter and $|x^m|_{t^w}^2=|t^{w_m}x^m|^2$, $|x|_{t^w}^2 = \sum_{m\in\Delta_5}|x^m|_{t^w}^2$.\\

With this moment map, it is convenient to assign $\mathbb{CP}^2$ with the associated toric \k metric $\omega_{t^w}$. The \k potential of $\omega_{t^w}$ is $\log |x|_{t^w}^2$.\\

$\Delta_5$ can also be viewed as a real triangle in $M_\mathbb{R} = M\otimes \mathbb{R}$. Then $w=(w_m)_{m\in \Delta_5}$ defines a function on the integral points in $\Delta_5$. If $w$ is convex, $w$ can be extended to a piecewise linear convex function on the real triangle $\Delta_5$. We will denote the extension also by $w$. A {\it generic} such $w$ will determine a simplicial decomposition of $\Delta_5$, with zero simplices being the integral points in $\Delta_5$. In this case, we say the piecewise linear convex function $w$ is compatible with the simplicial decomposition of $\Delta_5$, and the simplicial decomposition of $\Delta_5$ is determined by $w$. Conversely, as pointed out to me by Yi Hu, not every simplicial decompositions of $\Delta_5$ with zero simplices being integral points in $\Delta_5$ possesses a compatible piecewise linear convex function. We will only restrict our discussion to those simplicial decompositions that possess compatible piecewise linear convex functions. Consider the baricenter subdivision of such a simplicial decomposition of $\Delta_5$. Let $\Gamma_w$ denote the union of simplices in the baricenter subdivision that do not intersect integral points in $\Delta_5$. Then it is not hard to see that $\Gamma_w$ is a one-dimensional graph. $\Gamma_w$ divides $\Delta_5$ into regions, with a unique integral point of $\Delta_5$ located at the center of each of these regions. In particular we can think of the regions as parametrized by integral points in $\Delta_5$.\\

Let $C$ denote the quintic curve defined by the quintic polynomial

\[
p(x)= \sum_{|m|\leq 5}a_m x^m,
\]

with $|a_m|=t^{w_m}$. Combine theorem 3.1, proposition 5.7, theorem 5.1 and the remark following theorem 5.1 in \cite{N}, the theorem 9.1 in \cite{lag2} can be generalized in the general quintic case as the following.\\
\begin{theorem}
\label{bc}
For $w=(w_m)_{m\in \Delta_5}$ convex and positive, and $t$ small enough, $\tilde{\Gamma}_w = F_{t^w}(C)$ will be a fattening of a graph $\Gamma_w$. There exists a family $\{H_s\}_{s\in [0,1]}$ of piecewise smooth Lipschitz continuous symplectic automorphisms of $(\mathbb{CP}^2,\omega_{t^w})$ that restrict to identity on the three coordinate $\mathbb{CP}^1$'s, such that $H_0 = {\rm id}$, $\hat{C} = H_1(C)$ is a piecewise smooth symplectic curve satisfying $F_{t^w}(\hat{C})=\Gamma_w$.
\end{theorem}
\begin{flushright} $\Box$ \end{flushright}
{\bf Remark:} The moment map $F_{t^w}$ is invariant under the real 2-torus action. For any other moment map $F$ that is invariant under the real 2-torus action, $F(\hat{C})$ is a 1-dimensional graph.\\

Although not in the form we can directly apply, there are earlier results on the structures of amoeba with very different applications in mind that are more or less equivalent to theorem 3.1 of \cite{N}, which is stated in the first sentence of theorem \ref{bc}. References for these earlier works can be found in \cite{N} and \cite{M}.\\

{\bf Example:} For the standard simplicial decomposition of the Newton polygon $\Delta_5$ of quintic polyonmials (shown in Figure 1),\\
\begin{center}
\setlength{\unitlength}{1pt}
\begin{picture}(200,160)(-30,10)
\put(150,48){\line(-3,-5){18}}
\put(6,48){\line(3,-5){18}}
\put(60,138){\line(1,0){36}}
\put(132,78){\line(-3,-5){36}}
\put(24,78){\line(3,-5){36}}
\put(42,108){\line(1,0){72}}
\put(114,108){\line(-3,-5){54}}
\put(42,108){\line(3,-5){54}}
\put(24,78){\line(1,0){108}}
\put(96,138){\line(-3,-5){72}}
\put(60,138){\line(3,-5){72}}
\put(6,48){\line(1,0){144}}
\put(78,168){\line(-3,-5){90}}
\put(78,168){\line(3,-5){90}}
\put(-12,18){\line(1,0){180}}
\multiput(78,168)(36,0){1}{\circle*{4}}
\multiput(60,138)(36,0){2}{\circle*{4}}
\multiput(42,108)(36,0){3}{\circle*{4}}
\multiput(24,78)(36,0){4}{\circle*{4}}
\multiput(6,48)(36,0){5}{\circle*{4}}
\multiput(-12,18)(36,0){6}{\circle*{4}}
\end{picture}
\end{center}
\begin{center}
\stepcounter{figure}
Figure \thefigure: the standard simplicial decomposition
\end{center}
we have the corresponding $\hat{\Gamma} = F(C)$ (Figure 2),\\
\begin{center}
\begin{picture}(200,180)(-30,0)
\thicklines
\put(64,158){\line(3,-1){12}}
\put(76,154){\line(1,0){12}}
\put(100,158){\line(-3,-1){12}}

\multiput(0,0)(-18,-30){4}
{\put(64,158){\line(4,-3){15.2}}
\put(79.2,146.6){\line(0,-1){17}}
\put(46,128){\line(3,-1){18}}
\put(64,122){\line(2,1){15.4}}}

\multiput(164,0)(18,-30){4}
{\put(-64,158){\line(-4,-3){15.2}}
\put(-79.2,146.6){\line(0,-1){17}}
\put(-46,128){\line(-3,-1){18}}
\put(-64,122){\line(-2,1){15.4}}}

\put(-54,-90)
{\put(46,128){\line(4,-3){10}}
\put(56,120.5){\line(3,-5){6}}
\put(64,98){\line(-1,6){2.1}}}

\put(218,-90)
{\put(-46,128){\line(-4,-3){10}}
\put(-56,120.5){\line(-3,-5){6}}
\put(-64,98){\line(1,6){2.1}}}

\multiput(-54,-90)(36,0){4}
{\put(64,98){\line(1,6){3.3}}
\put(100,98){\line(-1,6){3.3}}
\put(82,125){\line(2,-1){15}}
\put(82,125){\line(-2,-1){15}}}

\multiput(82,108.5)(18,-30){3}{\circle{30}}
\multiput(64,78.5)(18,-30){2}{\circle{30}}
\multiput(46,48.5)(18,-30){1}{\circle{30}}

\thinlines
\multiput(82,149)(36,0){1}{\line(2,1){18}}
\multiput(82,149)(36,0){1}{\line(-2,1){18}}
\multiput(82,128)(36,0){1}{\line(0,1){21}}
\multiput(46,128)(36,0){2}{\line(2,-1){18}}
\multiput(82,128)(36,0){2}{\line(-2,-1){18}}
\multiput(64,98)(36,0){2}{\line(0,1){21}}
\multiput(28,98)(36,0){3}{\line(2,-1){18}}
\multiput(64,98)(36,0){3}{\line(-2,-1){18}}
\multiput(46,68)(36,0){3}{\line(0,1){21}}
\multiput(10,68)(36,0){4}{\line(2,-1){18}}
\multiput(46,68)(36,0){4}{\line(-2,-1){18}}
\multiput(28,38)(36,0){4}{\line(0,1){21}}
\multiput(-8,38)(36,0){5}{\line(2,-1){18}}
\multiput(28,38)(36,0){5}{\line(-2,-1){18}}
\multiput(10,8)(36,0){5}{\line(0,1){21}}
\put(-26,8){\line(1,0){216}}
\put(-26,8){\line(3,5){108}}
\put(190,8){\line(-3,5){108}}
\end{picture}
\end{center}
\begin{center}
\stepcounter{figure}
Figure \thefigure: $\tilde{\Gamma}$ for the standard simplicial decomposition
\end{center}
which is a fattening of the following graph $\Gamma$ (Figure 3). (We should remark here that Figure 2 is a rough topological illustration of the image $\hat{\Gamma} = F(C)$. Some part of the edges of the image that are straight or convex could be curved or concave in more accurate picture. Of course, such inaccuracy will not affect our mathematical argument and the fact that $\hat{\Gamma} = F(C)$ is a fattening of a graph.) By results in \cite{N}, we can symplectically deform $C$ to a symplectic curve $\hat{C}$ such that $F(\hat{C}) = \Gamma$.\\
\begin{center}
\setlength{\unitlength}{1pt}
\begin{picture}(200,160)(-30,10)
\put(150,48){\line(-3,-5){18}}
\put(6,48){\line(3,-5){18}}
\put(60,138){\line(1,0){36}}
\put(132,78){\line(-3,-5){36}}
\put(24,78){\line(3,-5){36}}
\put(42,108){\line(1,0){72}}
\put(114,108){\line(-3,-5){54}}
\put(42,108){\line(3,-5){54}}
\put(24,78){\line(1,0){108}}
\put(96,138){\line(-3,-5){72}}
\put(60,138){\line(3,-5){72}}
\put(6,48){\line(1,0){144}}
\put(78,168){\line(-3,-5){90}}
\put(78,168){\line(3,-5){90}}
\put(-12,18){\line(1,0){180}}
\multiput(78,168)(36,0){1}{\circle*{4}}
\multiput(60,138)(36,0){2}{\circle*{4}}
\multiput(42,108)(36,0){3}{\circle*{4}}
\multiput(24,78)(36,0){4}{\circle*{4}}
\multiput(6,48)(36,0){5}{\circle*{4}}
\multiput(-12,18)(36,0){6}{\circle*{4}}
\thicklines
\multiput(78,149)(36,0){1}{\line(2,1){18}}
\multiput(78,149)(36,0){1}{\line(-2,1){18}}
\multiput(78,128)(36,0){1}{\line(0,1){21}}
\multiput(42,128)(36,0){2}{\line(2,-1){18}}
\multiput(78,128)(36,0){2}{\line(-2,-1){18}}
\multiput(60,98)(36,0){2}{\line(0,1){21}}
\multiput(24,98)(36,0){3}{\line(2,-1){18}}
\multiput(60,98)(36,0){3}{\line(-2,-1){18}}
\multiput(42,68)(36,0){3}{\line(0,1){21}}
\multiput(6,68)(36,0){4}{\line(2,-1){18}}
\multiput(42,68)(36,0){4}{\line(-2,-1){18}}
\multiput(24,38)(36,0){4}{\line(0,1){21}}
\multiput(-12,38)(36,0){5}{\line(2,-1){18}}
\multiput(24,38)(36,0){5}{\line(-2,-1){18}}
\multiput(6,8)(36,0){5}{\line(0,1){21}}
\end{picture}
\end{center}
\begin{center}
\stepcounter{figure}
Figure \thefigure: $\Gamma$ for the standard simplicial decomposition
\end{center}
If the simplicial decomposition is changed to the following (Figure 4), we have the corresponding $\Gamma$ (Figure 4).\\
\begin{center}
\setlength{\unitlength}{1pt}
\begin{picture}(150,160)(20,10)
\put(61,136){\line(2,-1){52}}
\put(150,48){\line(-3,-5){18}}
\put(6,48){\line(3,-5){18}}
\put(60,138){\line(1,0){36}}
\put(132,78){\line(-3,-5){36}}
\put(24,78){\line(3,-5){36}}
\put(42,108){\line(1,0){72}}
\put(114,108){\line(-3,-5){54}}
\put(42,108){\line(3,-5){54}}
\put(24,78){\line(1,0){36}}
\put(96,78){\line(1,0){36}}
\put(78,48){\line(0,1){60}}
\put(78,108){\line(-3,-5){54}}
\put(60,138){\line(3,-5){72}}
\put(6,48){\line(1,0){144}}
\put(78,168){\line(-3,-5){90}}
\put(78,168){\line(3,-5){90}}
\put(-12,18){\line(1,0){180}}
\multiput(78,168)(36,0){1}{\circle*{4}}
\multiput(60,138)(36,0){2}{\circle*{4}}
\multiput(42,108)(36,0){3}{\circle*{4}}
\multiput(24,78)(36,0){4}{\circle*{4}}
\multiput(6,48)(36,0){5}{\circle*{4}}
\multiput(-12,18)(36,0){6}{\circle*{4}}
\end{picture}
\begin{picture}(150,160)(-20,10)
\put(61,136){\line(2,-1){52}}
\put(150,48){\line(-3,-5){18}}
\put(6,48){\line(3,-5){18}}
\put(60,138){\line(1,0){36}}
\put(132,78){\line(-3,-5){36}}
\put(24,78){\line(3,-5){36}}
\put(42,108){\line(1,0){72}}
\put(114,108){\line(-3,-5){54}}
\put(42,108){\line(3,-5){54}}
\put(24,78){\line(1,0){36}}
\put(96,78){\line(1,0){36}}
\put(78,48){\line(0,1){60}}
\put(78,108){\line(-3,-5){54}}
\put(60,138){\line(3,-5){72}}
\put(6,48){\line(1,0){144}}
\put(78,168){\line(-3,-5){90}}
\put(78,168){\line(3,-5){90}}
\put(-12,18){\line(1,0){180}}
\multiput(78,168)(36,0){1}{\circle*{4}}
\multiput(60,138)(36,0){2}{\circle*{4}}
\multiput(42,108)(36,0){3}{\circle*{4}}
\multiput(24,78)(36,0){4}{\circle*{4}}
\multiput(6,48)(36,0){5}{\circle*{4}}
\multiput(-12,18)(36,0){6}{\circle*{4}}
\thicklines
\multiput(78,149)(36,0){1}{\line(2,1){18}}
\multiput(78,149)(36,0){1}{\line(-2,1){18}}
\multiput(78,149)(4,-30){2}{\line(2,-3){14}}
\multiput(114,129)(-32,-10){2}{\line(-1,0){23}}
\multiput(82,119)(-27,-10){1}{\line(1,1){10}}
\multiput(42,128)(36,0){1}{\line(2,-1){18}}
\multiput(78,128)(36,0){0}{\line(-2,-1){18}}
\multiput(60,98)(36,0){1}{\line(0,1){21}}
\multiput(24,98)(72,0){2}{\line(2,-1){18}}
\multiput(60,98)(26,-20){2}{\line(1,-2){10}}
\put(70,78){\line(1,0){16}}
\multiput(60,98)(72,0){2}{\line(-2,-1){18}}
\multiput(96,98)(-26,-20){2}{\line(-1,-2){10}}
\multiput(42,68)(72,0){2}{\line(0,1){21}}
\multiput(6,68)(36,0){2}{\line(2,-1){18}}
\multiput(114,68)(36,0){1}{\line(2,-1){18}}
\multiput(42,68)(36,0){1}{\line(-2,-1){18}}
\multiput(114,68)(36,0){2}{\line(-2,-1){18}}
\multiput(24,38)(36,0){4}{\line(0,1){21}}
\multiput(-12,38)(36,0){5}{\line(2,-1){18}}
\multiput(24,38)(36,0){5}{\line(-2,-1){18}}
\multiput(6,8)(36,0){5}{\line(0,1){21}}
\end{picture}
\end{center}
\begin{center}
\stepcounter{figure}
Figure \thefigure: alternative simplicial decomposition and corresponding $\Gamma$
\end{center}
With this preparation, now we are ready to address the meaning of {\bf near the large complex limit}. We will resume our usual notation instead of the 2-dimensional notation at this point. Consider the quintic Calabi-Yau family $\{X_\psi\}$ in $\mathbb{CP}^4$ defined by\\
\[
p_\psi(z) = p_a(z) + \psi \prod_{k=1}^5 z_k = \sum_{m\not=m_0,|m|=5}a_mz^m + \psi \prod_{k=1}^5 z_k=0,
\]\\
where $|a_m| = t^{w_m}$, $t>0$ is a parameter and $\{w_m\}_{m\in \Delta}$ is a function defined on the Newton polyhedron $\Delta$ of the quintic polynomials on $\mathbb{CP}^4$, which is\\
\begin{equation}
\label{bh}
\Delta =\{z^m = z_1^{m_1}\cdots z_5^{m_5} | |m|=5,\ m=(m_1,\cdots,m_5) \in\mathbb{Z}_{\geq 0}^5\}.
\end{equation}
$\Delta$ is a reflexive polyhedron in the lattice\\
\begin{eqnarray}
\label{bi}
M &=& \{ m=(m_1,\cdots,m_5) \in\mathbb{Z}^5| |m|=5\}\\\nonumber &\cong& \{ m=(m_1,\cdots,m_5) \in\mathbb{Z}^5| |m|=0\}
\end{eqnarray}
identifying $m_0=(1,1,1,1,1)$ with the origin.\\

Let $\Delta^0$ denote the set of integral points in the 2-skeleton of $\Delta$. Or in other words, the integral points in $\Delta$ that are not in the interior of 3-faces and not the center $m_0=(1,1,1,1,1)$. For $m\in \Delta$, let $w'_m = w_m-w_{m_0}$ (recall $\psi = a_{m_0}$). We have two functions $w=(w_m)_{m\in \Delta}$ and $w'=(w'_m)_{m\in \Delta}$ defined on the integral points in $\Delta$.\\
\begin{de}
\label{bb}
$w'=(w'_m)_{m\in \Delta}$ is called convex with respect to  $\Delta^0$, if for any $m' \in \Delta^0$ there exists a linear function $n$ such that $n(m')= w'_{m'}$ and $n(m)\leq w'_m$ for $m\in \Delta^0\backslash \{m'\}$.
\end{de}
The two convexities we defined are closely related. For $I \subset \{1,2,3,4,5\}$ with $|I|=2$, let

\[
\Delta_I =\{z^{m_I} = \prod_{i\not\in I}z_i^{m_i} | |m_I|=5,\ m_I=(m_i)_{i\not\in I} \in\mathbb{Z}_{\geq 0}^3\}.
\]

Then

\[
\Delta^0 = \bigcup_{\tiny{\begin{array}{c}I \subset \{1,2,3,4,5\}\\ |I|=2\end{array}}} \Delta_I.
\]

Let

\[
w_I = w|_{\Delta_I}.
\]

Then we have the following lemma

\begin{lm}
$w'=(w'_m)_{m\in \Delta}$ is convex with respect to  $\Delta^0$ for $-w_{m_0}$ large enough, if and only if $w_I$ is convex on $\Delta_I$ in the sense of definition \ref{ba} for all $I \subset \{1,2,3,4,5\}$ with $|I|=2$.\\
\end{lm}
{\bf Proof:} The ``only if'' part is trivial. We will prove the ``if'' part. For any $m' \in \Delta^0$, there exists the smallest subface $\Delta_{m'} \subset \Delta$, such that $m' \in \Delta_{m'}$. Since $\Delta$ is convex, there exists a linear function $n_{m'}$ such that\\
\[
n_{m'}|_{\Delta_{m'}}=-1,\ \ n_{m'}|_{\Delta\backslash\Delta_{m'}} >-1.
\]\\
By assumption, $w_I$ is convex on $\Delta_I$ in the sense of definition \ref{ba} for all  $I \subset \{1,2,3,4,5\}$ with $|I|=2$. Clearly, $\Delta_{m'}$ belongs to one of $\Delta_I$. Therefore $w|_{\Delta_{m'}}$ is convex on $\Delta_{m'}$ in the sense of definition \ref{ba}. Namely, there exists a linear function $n'$ such that\\
\[
n'(m') = w_{m'},\ {\rm and}\ n'(m') \leq w_{m}\ \ {\rm for} \ m\in \Delta_{m'}\backslash \{m'\}.
\]\\
Let $n = n' + w_{m_0}n_{m'}$, then it is easy to see that for $-w_{m_0}$ large enough, $n(m')= w'_{m'}$ and $n(m)\leq w'_m$ for $m\in \Delta^0\backslash \{m'\}$, namely, $w'=(w'_m)_{m\in \Delta}$ is convex with respect to  $\Delta^0$.
\begin{flushright} $\Box$ \end{flushright}
With this lemma in mind, we can make the concept {\bf near the large complex limit} more precise in the following definition.\\
\begin{de}
\label{bd}
The quintic Calabi-Yau hypersurface $X_\psi$ is said to be {\bf near the large complex limit}, if $w'=(w'_m)_{m\in \Delta}$ is convex with respect to  $\Delta^0$ and $t$ is small.\\
\end{de}
When $X_\psi$ is near the large complex limit, for $I \subset \{1,2,3,4,5\}$ with $|I|=2$, the corresponding $w_I$ is convex on $\Delta_I$. When $w$ is generic, by previous construction, $w_I$ determines a 1-dimensional graph $\Gamma_I$ on $\Delta_I$. Let

\[
\Gamma = \bigcup_{\tiny{\begin{array}{c}I \subset \{1,2,3,4,5\}\\ |I|=2\end{array}}} \Gamma_I.
\]

At this point, in order to apply theorem \ref{bc}, we need to adopt the toric metric with the \k potential $\log |x|_{t^w}^2$ on $\mathbb{CP}^4$. The corresponding moment map is\\
\[
F_{t^w}(x) = \sum_{m\in\Delta^0} \frac{|x^m|_{t^w}^2}{|x|_{t^w}^2}m,
\]\\
which maps $\mathbb{CP}^4$ to $\Delta$, where $|x^m|_{t^w}^2=|t^{w_m}x^m|^2$, $|x|_{t^w}^2 = \sum_{m\in\Delta^0}|x^m|_{t^w}^2$.\\

By theorem \ref{bc}, $\tilde{\Gamma}_I = F_{t^w}(C_I)$ is a fattening of $\Gamma_I$. Therefore,
\[
\tilde{\Gamma} = \bigcup_{\tiny{\begin{array}{c}I \subset \{1,2,3,4,5\}\\ |I|=2\end{array}}} \tilde{\Gamma}_I
\]\\
is a fattening of 1-dimensional graph $\Gamma$. There are 3 types of singular fibers over different part of $\tilde{\Gamma}$. Let $\tilde{\Gamma}^2$ denote the interior of $\tilde{\Gamma}$, and\\
\[
\tilde{\Gamma}^1 = \partial \tilde{\Gamma} \cap \left(\bigcup_{|I|=2}\Delta_I\right),
\ \ \ \
\tilde{\Gamma}^0 = \partial \tilde{\Gamma} \cap \left(\bigcup_{|I|=3}\Delta_I\right).
\]\\
Then

\begin{equation}
\label{bg}
\tilde{\Gamma} = \tilde{\Gamma}^0 \cup \tilde{\Gamma}^1 \cup \tilde{\Gamma}^2.
\end{equation}

Let $F_\infty = F_{t^w}|_{X_\infty}: X_\infty \rightarrow \partial \Delta$ be the moment map Lagrangian torus fibration. Results in \cite{M,N} imply that the inverse image of $F_\infty$ over a point in $\tilde{\Gamma}^0$ is a circle that intersects $C$ at one point, over a point in $\tilde{\Gamma}^1$ is a 2-torus that intersects $C$ at one point, over a point in $\tilde{\Gamma}^2$ is a 2-torus that intersects $C$ at two points. Then by theorem \ref{bf}, we have the following theorem.\\
\begin{theorem}
When $X_\psi$ is near the large complex limit, starting with the topologically smooth Lagrangian torus fibration $F_\infty: X_\infty \rightarrow \partial \Delta$, the gradient flow method will produce a Lagrangian torus fibration $F_\psi: X_{\psi} \rightarrow \partial\Delta$ with codimension 1 singular locus $\tilde{\Gamma} = F_\infty(C)$ as a fattening of the graph $\Gamma$. There are 4 types of fibers.\\
(i). For $p\in \partial\Delta\backslash \tilde{\Gamma}$, $F_\psi^{-1}(p)$ is a smooth Lagrangian 3-torus.\\
(ii). For $p\in \tilde{\Gamma}^2$, $F_\psi^{-1}(p)$ is a Lagrangian 3-torus with two circles collapsed to two singular points.\\
(iii). For $p\in \tilde{\Gamma}^1$, $F_\psi^{-1}(p)$ is a Lagrangian 3-torus with one circle collapsed to one singular point.\\
(iv). For $p\in \tilde{\Gamma}^0$, $F_\psi^{-1}(p)$ is a Lagrangian 3-torus with one 2-torus collapsed to one singular points.
\end{theorem}
\begin{flushright} $\Box$ \end{flushright}
Let $\Gamma = \Gamma^1\cup \Gamma^2\cup \Gamma^3$, where $\Gamma^1$ is the smooth part of $\Gamma$, $\Gamma^2$ is the singular part of $\Gamma$ in the interior of the 2-skeleton of $\Delta$, $\Gamma^3$ is the singular part of $\Gamma$ in the 1-skeleton of $\Delta$. Applying theorem \ref{bf}, with help of the second part of the theorem \ref{bc} and some other results in \cite{lag2,N}, we can construct a Lagrangian torus fibration $\hat{F}_{\psi}: X_{\psi} \rightarrow \partial\Delta$ with codimension 2 singular locus $\Gamma = \hat{F}_\psi(C)$ in the following theorem. For the definition of type I, II, III singular fibers, please refer to the section 8.\\
\begin{theorem}
\label{be}
When $X_\psi$ is near the large complex limit, the gradient flow method will produce a Lagrangian torus fibration $\hat{F}_{\psi}: X_{\psi} \rightarrow \partial\Delta$ with the graph singular locus $\Gamma = \hat{F}_\psi(C)$. There are 4 types of fibers.\\
(i). For $p\in \partial\Delta\backslash \Gamma$, $\hat{F}_{\psi}^{-1}(p)$ is a smooth Lagrangian 3-torus.\\
(ii). For $p\in \Gamma^1$, $\hat{F}_{\psi}^{-1}(p)$ is a type I singular fiber.\\
(iii). For $p\in \Gamma^2$, $\hat{F}_{\psi}^{-1}(p)$ is a type II singular fiber.\\
(iv). For $p\in \Gamma^3$, $\hat{F}_{\psi}^{-1}(p)$ is a type III singular fiber.\\
\end{theorem}
{\bf Proof:} According to theorem \ref{bc}, on each $\overline{D_I} \cong \mathbb{CP}^2$ in the 2-skeleton of $X_\infty$, there exists a family $\{H_s\}_{s\in [0,1]}$ of piecewise smooth Lipschitz continuous symplectic automorphisms of $\mathbb{CP}^2$ that restrict to identity on the three coordinate $\mathbb{CP}^1$'s, such that $H_0 = {\rm id}$, $\hat{C}_I = H_1(C_I)$ is a piecewise smooth symplectic curve satisfying $F_\infty(\hat{C}_I)=\Gamma_I$. These $\{H_s\}_{s\in [0,1]}$'s can be patched together to form a family of symplectic automorphisms of the 2-skeleton of $X_\infty$. Use the extension theorem (theorem 6.8) from \cite{lag2}, we may extend such automorphisms to $X_\infty$ to form a family (still denote by $\{H_s\}_{s\in [0,1]}$) of piecewise smooth Lipschitz continuous symplectic automorphisms of $X_\infty$ such that $H_0 = {\rm id}$, $\hat{C} = H_1(C)$ is a normal crossing union of piecewise smooth symplectic curves satisfying $F_\infty(\hat{C})=\Gamma$.\\

Clearly $\hat{F}_\infty = F_\infty \circ H_1: X_\infty \rightarrow \partial \Delta$ is also a topologically smooth Lagrangian fibration that satisfies $\hat{F}_\infty(C) = \Gamma$. Apply theorem \ref{bf} to the topologically smooth Lagrangian fibration $\hat{F}_\infty: X_\infty \rightarrow \partial \Delta$, we will get the desired result.
\begin{flushright} $\Box$ \end{flushright}
{\bf Remark:} Although the Lagrangian fibrations produced by this theorem have codimension 2 singular locus, this fibration map is not $C^{\infty}$ map (merely Lipschitz). I believe that a small perturbation of this map will give an almost $C^{\infty}$ Lagrangian fibration with the same topological structure. (At least $C^{\infty}$ away from type II singular fibers.) In \cite{smooth}, we will discuss partial smoothing of the Lagrangian torus fibration constructed in theorem \ref{be}.\\

\se{Complex moduli, \k moduli and monomial-divisor mirror map}
In this section we discuss the partial compactification and chamber decomposition of the complex moduli space of quintics near the large complex limit based on the structure of Lagrangian torus fibrations. We will show that this partial compactification is equivalent to the well known secondary fan compactification (\cite{AGM,AGM2}). We will also review the quotient construction for the mirror of quintics, the \k moduli of the mirror of quintics and the monomial-divisor mirror map. These will be used in the later section for establishing the symplectic SYZ correspondence. Most materials in this section are well known (see \cite{GP, Roan, AGM, AGM2, mirrorbook}). The main purpose of our presentation here is for reviewing the facts and fixing the notations, with the exception of our new interpretation of secondary fan compactification in terms of Lagrangian torus fibration and the slicing theorem, which removes some restriction in the results of \cite{AGM} so that they can be applied in our situation.\\

\subsection{Slicing theorem and the complex moduli space of quintic Calabi-Yau manifolds}
For quintic Calabi-Yau hypersurfaces, the complex moduli can be seen as the space of homogeneous quintic polynomials on $\mathbb{C}^5$ modulo the action of $GL(5,\mathbb{C})$. Recall from (\ref{bh}) and (\ref{bi}) the Newton polyhedron $\Delta \subset M$ of quintic polynomials. The space of quintic polynomials $Q$ can be viewed as\\
\[
Q ={\rm Span}_\mathbb{C}(\Delta).
\]\\
To understand the quotient of $Q$ by $GL(5,\mathbb{C})$, classically we would have to use the geometric invariant theory. It is fairly complicated, and the compactification is not unique.
For the purpose of mirror symmetry, the GIT approach is a bit unnatural. For example, the large complex limit

\[
p_\infty(z) =  \prod_{k=1}^5 z_k
\]

is not invariant under the action of $GL(5,\mathbb{C})$. The stablizer of $p_\infty(z)$ (at least infinitesimally) is merely the Cartan subgroup $T\in SL(5,\mathbb{C})$. On the other hand, mirror symmetry mainly concerns the behavior of the complex moduli near the large complex limit $p_\infty(z)$, with the meaning of ``near" suitably specified.\\

For these reasons, it will be much nicer if we can find a canonical $T$-invariant slice $Q_0$ in $Q$ that passes through the large complex limit $p_\infty(z)$ and intersects each orbit of $GL(5,\mathbb{C})$ exactly at an orbit of $T$. This seems quite impossible to be done on the whole complex moduli. But near the large complex limit $p_\infty(z)$, this can be done very nicely by the following slicing theorem.\\

Recall that $\Delta^0$ is the two skeleton of $\Delta$. Let

\[
Q_0 = m_0 + {\rm Span}_\mathbb{C}(\Delta^0).
\]

Then we have

\begin{theorem}
Near the large complex limit $p_\infty$, $Q_0 = m_0 + {\rm Span}_\mathbb{C}(\Delta^0)$ is a slice of $Q ={\rm Span}_\mathbb{C}(\Delta)$ that contains the large complex limit $p_\infty$, and is invariant under toric automorphism group $T$ and intersects each nearby orbit of the automorphism group $GL(5,\mathbb{C})$ at a finite set of orbits of the toric automorphism group $T$ (the maximal torus in $SL(5,\mathbb{C})$) that is parametrized by a quotient of the Weyl group of $SL(5,\mathbb{C})$.\\
\end{theorem}
{\bf Proof:} Let\\
\[
F(z)= \sum_{m(\not=m_0)\in \Delta} a_mz^m + \psi z^{m_0},
\]\\
then the theorem asserts that when $\psi$ is large, through linear transformation on $z$, $F(z)$ can be reduced to the following standard form\\
\[
F_0(z)=\sum_{m\in \Delta^0} a_mz^m + \psi z^{m_0}.
\]\\
In general, it seems a rather hard question. But when $\psi$ is large (when near the large complex limit), the question is much easier. Let us consider the linear transformation\\
\[
L_1:z\rightarrow (I- \frac{1}{\psi}B)z,
\]\\
where $b_{jk} = a_{m_{jk}}$, $m_{jk}=m_0-e_j+e_k$. Then\\
\begin{eqnarray*}
F_1(z)=F(L_1z) &=& \psi z^{m_0} + \sum_{m(\not=m_0)\in \Delta} a_mz^m - \sum_{j,k}b_{jk}z^{m_{jk}} + O(\frac{1}{\psi})\\
&=& \psi z^{m_0} + \sum_{m\in \Delta^0} a_mz^m + \sum_{j,k} a_{m_{jk}}z^{m_{jk}} - \sum_{j,k}b_{jk}z^{m_{jk}} + O(\frac{1}{\psi})\\
&=& \psi z^{m_0} + \sum_{m\in \Delta^0} a_mz^m + O(\frac{1}{\psi})
\end{eqnarray*}
Adjust $\psi$ and $a_m$ accordingly, above equation can also be expressed as\\
\[
F_1(z)= \sum_{m(\not=m_0)\in \Delta} a_mz^m + \psi z^{m_0},
\]\\
where $a_m=O(1/\psi)$ for $m\in \Delta\backslash \Delta^0$. Repeating above process, consider\\
\[
L_2:z\rightarrow (I- \frac{1}{\psi}B)z.
\]\\
Then\\
\[
F_2(z) = F_1(L_2z)=F(L_1L_2z)= \sum_{m(\not=m_0)\in \Delta} a_mz^m + \psi z^{m_0},
\]\\
where $a_m=O(1/\psi^2)$ for $m\in \Delta\backslash \Delta^0$. Repeat this process inductively, and define\\
\[
L= L_1L_2\cdots.
\]
Since $L_l-I = O(1/\psi^l)$, the infinite product converges when $\psi$ is large. Then\\
\[
F_0(z) = F(Lz) = \sum_{m\in \Delta^0} a_mz^m + \psi z^{m_0}
\]\\
is in standard form (belonging to our slice after divided by $\psi$).\\

To show that near $p_\infty$ $Q_0$ intersects each $GL(5,\mathbb{C})$ orbit at an orbit of $T$, notice that any group element in $GL(5,\mathbb{C})\backslash T$ near $T$, when acting on elements of $Q_0$ near $p_\infty$, will produce terms in $\Delta\backslash (\Delta^0\cup \{m_0\})$ or alter the $m_0$ term. This fact together with the fact that stablizer of $p_\infty$ is exactly the normalizer of $T\subset SL(5,\mathbb{C})$ imply our conclusion.
\begin{flushright} $\Box$ \end{flushright}
{\bf Remark:} The similar slicing theorem is also true in the general toric case, which will be discussed in \cite{tor}. Given the elementary nature of the theorem, we believe it may have appeared in the literature in some form, although we are not aware of it. In any case, this kind of slicing theorem is very important in understanding the complex moduli of Calabi-Yau hypersurfaces in a toric variety near the large complex limit.\\

As we know, the normalizer of $T$ in $SL(5,\mathbb{C})$ is a semi-direct product of $T$ and the Weyl group $S_5$ of $SL(5,\mathbb{C})$. The construction of the complex moduli near the large complex limit involves first constructing certain ``quotient" of $Q_0$ by $T$ that is invariant under $S_5$, then taking the further quotient by $S_5$.\\

$Q_0$ is naturally an affine space with the natural toric structure determined by the open complex torus $T_{Q_0} \subset Q_0$. The action of $T$ is compatible with this toric structure. We expect the ``quotient" of $Q_0$ by $T$ to be a toric variety as toric partial compactification of the quotient complex torus $T_{Q_0}/T$. Before getting into the detail, we first introduce some notations. Let\\
\[
\tilde{M}_0 = \left\{a^I= \prod_{m\in \Delta^0}a_m^{i_m}\ \left|\ I=(i_m)_{m\in \Delta^0}\in \mathbb{Z}^{\Delta^0}\right.\right\}\cong \mathbb{Z}^{\Delta^0},
\]\\
then its dual $\tilde{N}_0$ is naturally the space of all the weights\\
\[
\tilde{N}_0 = \{ w=(w_m)_{m\in \Delta^0}\in \mathbb{Z}^{\Delta^0} \} \cong \mathbb{Z}^{\Delta^0}.
\]\\
Follow the toric geometry convention, let $N$ denote the dual lattice of $M$. $n\rightarrow w_m = \langle m,n\rangle$ defines an embedding $N\hookrightarrow \tilde{N}_0$. Let $W\in N_0$ be the image of this embedding, $\tilde{M}=W^\perp$, then $\tilde{N} = \tilde{N}_0/W$ is dual to $\tilde{M}$. $\tilde{M}$ can be viewed as the lattice of monomials on the quotient complex torus $T_{Q_0}/T$. Therefore, the fan of the ``quotient" of $Q_0$ by $T$ should be in $\tilde{N}$. For a convex cone $\sigma_0\subset \tilde{N}_0$, let $\sigma^\vee = \sigma_0^\vee\cap \tilde{M}$. Then $\sigma = (\sigma^\vee)^\vee$ is naturally the projection of $\sigma_0$ to $\tilde{N}$.\\

Our discussion on the notion of near the large complex limit (definition \ref{bd}), which was needed for the construction of the Lagrangian fibration, naturally suggests a chamber decomposition of the complex moduli and the related fan structure in $\tilde{N}$ that determines a ``canonical'' partial compactification of the complex moduli near the large complex limit.\\

Recall from definition \ref{bd}, the quintic hypersurface $X$ defined by

\[
p(z) = z^{m_0} + \sum_{m\in \Delta^0} a_m z^m \in T_{Q_0}
\]

with $|a_m| = t^{w_m}$ is near the large complex limit if $\{w_m\}_{m\in \Delta^0}$ is strictly convex and $t>0$ is small. (A generic convex $w=(w_m)_{m\in \Delta^0}$ will determine a simplicial decomposition $Z$ of $\Delta^0$ with integral points in $\Delta^0$ as vertices of the simplices in $Z$. In such situation, we say $w=\{w_m\}_{m\in \Delta^0}$ is compatible with $Z$.) Notice that in our discussion, we are repeatedly using the notation $\Delta^0$ for two slightly different meanings, first, as the set of integral points in the two skeleton of polyhedron $\Delta$, secondly, as the real two skeleton of polyhedron $\Delta$.\\

According to theorem \ref{be}, the gradient flow method will produce a Lagrangian torus fibration for $X$ with singular locus $\Gamma_Z$, which is determined by the simplicial decomposition $Z$ of $\Delta^0$ that is compatible with the strictly convex function $\{w_m\}_{m\in \Delta^0}$. For each $Z$, those quintic polynomials in $T_{Q_0}$ determining the same $Z$ form a chamber $\tilde{U}_Z$ in $T_{Q_0}$. (The closure of $\tilde{U}_Z$ forms a chamber in $Q_0$.) For Calabi-Yau quintics in $\tilde{U}_Z$ near the large complex limit, the corresponding Lagragian fibrations have singular locus with the same combinatorial type $\Gamma_Z$. Since the convex function $\{w_m\}_{m\in \Delta^0}$ modulo linear functions on $M$ determine the same $Z$, $\tilde{U}_Z$ is invariant under the $T$ action. $U_Z = \tilde{U}_Z/T$ forms a chamber in the quotient complex torus $T_{Q_0}/T$.\\

Let $\tilde{Z}$ denote the set of all simplicial decomposition $Z$ of $\Delta^0$ with integral points in $\Delta^0$ as vertices of the simplices in $Z$. Then we get a chamber decomposition $\{\tilde{U}_Z\}_{Z\in\tilde{Z}}$ of $T_{Q_0}$ near the large complex limit $m_0$ that is invariant under $T$ and permuted by the action of the Weyl group of $SL(5,\mathbb{C})$, hence decends to a chamber decomposition $\{U_Z\}_{Z\in\tilde{Z}}$ of the quotient complex torus $T_{Q_0}/T$, whose quotient by the Weyl group of $SL(5,\mathbb{C})$ form a Zariski open set of the the complex moduli of the quintic.\\

On the other hand, for any $Z \in \tilde{Z}$,

\[
\tau_Z = \{ w=(w_m)_{m\in \Delta^0}\in \tilde{N} | w\ {\rm is\ convex\ and\ compatible\ with\ } Z \}/W
\]

is a top dimensional cone in $\tilde{N}$. Let $\tilde{\Sigma}$ denote the fan consisting of subcones of $\tau_Z$ for $Z\in \tilde{Z}$. Then the support of the fan $\tilde{\Sigma}$ is\\
\[
\tau = \bigcup_{Z\in \tilde{Z}} \tau_Z.
\]\\
The corresponding toric variety $P_{\tilde{\Sigma}}$ is our model for the ``quotient" of $Q_0$ by $T$ as partial compactification of the quotient complex torus $T_{Q_0}/T$, which is compatible with the chamber decomposition in the sense that $U_Z$ can be naturally identified with the complexification of $\tau_Z$. It is easy to see that $\tilde{\Sigma}$ is invariant under the natural action of the Weyl group of $SL(5,\mathbb{C})$. The actual complex moduli space of quintic near the large complex limit is $P_{\tilde{\Sigma}}$ modulo the action of the Weyl group $S_5$ of $SL(5,\mathbb{C})$.\\
\begin{theorem}
$P_{\tilde{\Sigma}}/S_5$ is a natural compactification of the moduli space of quintic Calabi-Yau hypersurfaces near the large complex limit.
\end{theorem}
\begin{flushright} $\Box$ \end{flushright}
{\bf Remark:} $\tilde{\Sigma}$ is the so-called secondary fan. Our discussion in this subsection provides an intrinsic geometric reconstruction of the secondary fan compactification based on the construction of Lagrangian torus fibration in the case of quintics. Other more detailed justification for the general case of Calabi-Yau hyperserfaces in toric variety will be discussed in \cite{tor}.\\

\subsection{The mirror of quintics and their \k moduli}
All the materials in this subsection are well known and have appeared in the literature (\cite{Roan,GP,Can,AGM,AGM2, Bat3}). For reader's convenience and fixing the notations for our application, we will review these results in the form we need in this subsection. We first recall the quotient construction for the mirror of quintics.\\

Let $(e_1,\cdots,e_5)$ be the standard basis of $M_0\cong \mathbb{Z}^{5}$, and $(e^1,\cdots,e^5)$ be the dual basis of $N_0\cong \mathbb{Z}^{5}$. Recall from (\ref{bh}) and (\ref{bi}) the Newton polyhedron $\Delta \subset M$ of quintic polynomials.\\
\[
N = M^\vee = \{ n=(n_1,\cdots,n_5) \in\mathbb{Z}^5\}/\mathbb{Z}n_0
\]\\
with $n_0=(1,1,1,1,1)$ is naturally the dual lattice of $M$. Let $\Delta^\vee$ be the dual polyhedron of $\Delta$ in $N$. Then the simplex $\Delta^\vee\subset N$ has vertices $\{n^i\}_{i=1}^5$, and the dual simplex $\Delta \subset M$ has vertices $\{m^i\}_{i=1}^5$, where\\
\[
n^i=[e^i],\ \ \ m^i=5e_i- m_0.
\]\\
Interestingly, $\{n^i\}_{i=1}^5$ and $\{m^i\}_{i=1}^5$ satisfy the same linear relation:\\
\[
\sum_{k=1}^5 n^k=0,\ \ \ \sum_{k=1}^5 m^k=0.
\]\\
Since $N$ is generated by $\{n^i\}_{i=1}^5$, there is a natural map $Q: N\rightarrow M$ that maps $n^i$ to $m^i$. The following lemma is very straightforward to check.\\
\begin{lm}
\[
M/Q(N) \cong (\mathbb{Z}_5)^3
\]
\end{lm}
\begin{flushright} $\Box$ \end{flushright}
Let $\Sigma_\Delta$ ($\Sigma_{\Delta^\vee}$) be the fan in $N$ ($M$) whose cones are spanned by faces of $\Delta^\vee$ ($\Delta$) from the origin. The corresponding toric variety $P_{\Sigma_\Delta}$ ($P_{\Sigma_{\Delta^\vee}}$) is the anti-canonical model, namely the anti-canonical class is ample on $P_{\Sigma_\Delta}$ ($P_{\Sigma_{\Delta^\vee}}$). The map $Q: N\rightarrow M$ gives us the well known equivalence between the Batyrev mirror construction (\cite{Bat3}) and the quotient interpretation of the mirror (\cite{Roan,GP,Can}).\\
\begin{co}
\[
P_{\Sigma_{\Delta^\vee}} \cong P_{\Sigma_\Delta}/(\mathbb{Z}_5)^3
\]
\end{co}
\begin{flushright} $\Box$ \end{flushright}
In this way, Calabi-Yau hypersurfaces in $P_{\Sigma_{\Delta^\vee}}$ via pullback are equivalent to $(\mathbb{Z}_5)^3$-invariant Calabi-Yau hypersurfaces in $P_{\Sigma_\Delta}$, namely, the Fermat type quintics $X_\psi$ defined by\\
\[
\sum_{k=1}^5 z^5_k +5\psi \prod_{k=1}^5z_k =0.
\]\\
Denote the quotient of $X_\psi$ in $P_{\Sigma_{\Delta^\vee}}$ to be $Y_\psi^0$. Then $Y_\psi^0$ corresponds to the mirror of quintics, equivalently according to both the Batyrev and the quotient mirror constructions. We are interested in Calabi-Yau hypersurfaces near the large complex limit, which correspond to large $\psi$.\\

Calabi-Yau hypersurfaces in $P_{\Sigma_{\Delta^\vee}}$ are not exactly the mirrors, since they are singular. To get smooth mirror Calabi-Yau manifolds, it is necessary to do some crepant resolution. This can be done by resolving the singularities of $P_{\Sigma_{\Delta^\vee}}$ suitably, then the pullback of $Y_\psi^0$ to the resolution, denoted by $Y_\psi$, gives the corresponding smooth Calabi-Yau manifold.\\

In general resolution of singularities of $P_{\Sigma_{\Delta^\vee}}$ corresponds to subdivision of the fan $\Sigma_{\Delta^\vee}$. Arbitrary resolution of $P_{\Sigma_{\Delta^\vee}}$ can destroy the Calabi-Yau property of $Y_\psi$. Those admissible resolutions, which preserve the Calabi-Yau property of $Y_\psi$, correspond to those subdivision fan $\Sigma^\vee$ of $\Sigma_{\Delta^\vee}$ whose generators of 1-dimensional cones lie in $\partial \Delta$. In other words, admissible resolutions of $P_{\Sigma_{\Delta^\vee}}$ correspond to polyhedron subdivision of $\partial \Delta$ with vertices of the polyhedrons being the integral points in $\Delta$. Different subdivisions $\Sigma^\vee$ of $\Sigma_{\Delta^\vee}$ correspond to different birational models $P_{\Sigma^\vee}$ of $P_{\Sigma_{\Delta^\vee}}$.\\

For each particular birational model $P_{\Sigma^\vee}$ with the birational map $\hat{\pi}: P_{\Sigma^\vee} \rightarrow P_{\Sigma_{\Delta^\vee}}$, the K\"{a}hler cone of $P_{\Sigma^\vee}$ can be understood as the set of piecewise linear convex functions on $M$ that are compatible with $\Sigma^\vee$ modulo linear functions on $M$. We are actually interested in the K\"{a}hler cone of the corresponding Calabi-Yau hypersurface $Y_\psi = \hat{\pi}^{-1}(Y_\psi^0) \subset P_{\Sigma^\vee}$. By K\"{a}hler cone of $Y_\psi$ in this paper, we only refer to those K\"{a}hler forms of $Y_\psi$ that come from the restriction of K\"{a}hler forms on $P_{\Sigma^\vee}$, although presumablly, there might be other K\"{a}hler forms not comming from such restriction.\\

As a toric variety, $P_{\Sigma_{\Delta^\vee}}$ can be seen as a union of complex tori of various dimensions (orbits of $T_M$). Under the moment map, these complex tori are mapped to the interior of subfaces of $\Delta^\vee$ of various dimensions. There is a dual relation from subfaces of $\Delta$ to subfaces of $\Delta^\vee$. For $\alpha$ a subface of $\Delta$, the dual face is\\
\[
\alpha^* = \{n\in \Delta^\vee |\langle m,n\rangle =-1\}.
\]\\
Clearly, $(\alpha^*)^* = \alpha$ and $\dim \alpha +\dim \alpha^* =3$. $P_{\Sigma^\vee}$ is a toric resolution of $P_{\Sigma_{\Delta^\vee}}$. Tori divisors in $P_{\Sigma^\vee}$ that dominate the complex torus in $P_{\Sigma_{\Delta^\vee}}$ corresponding to subface $\alpha$ in $\Delta^\vee$ are parametrized by interior integral points in $\alpha^*$.\\

Notice that $Y_\psi^0$ does not intersect the points (0-dimensional tori) in $P_{\Sigma_{\Delta^\vee}}$ corresponding to vertices (0-dimensional subfaces) of $\Delta^\vee$. A toric divisor in $P_{\Sigma^\vee}$ restricts non-trivially to $Y_\psi$ if and only if the corresponding integral point in $\Delta$ is in the 2-skeleton $\Delta^0 \subset \Delta$. Let $Z$ denote the simplicial decomposition of $\Delta^0$ determined by the fan $\Sigma^\vee$. Then we have\\
\begin{prop}
\label{ca}
The K\"{a}hler cone of $Y_\psi \subset P_{\Sigma^\vee}$ can be naturally identified as the set of $w=(w_m)_{m\in \Delta^0}$ (understood as a piecewise linear function on $M$) that is convex with respect to $\Delta^0$ in the sense of definition \ref{bd} and compatible with $Z$, modulo linear functions on $M$.
\end{prop}
\begin{flushright} $\Box$ \end{flushright}
Notice that $Y_\psi$ is the same for all $\Sigma^\vee$ that determine the same simplicial decomposition $Z$ of $\Delta^0$, and the K\"{a}hler cone of $Y_\psi$ determined in proposition \ref{ca} is exactly the $\tau_Z$ as introduced in the previous subsection. From the point of view of mirror symmetry, it is more natural to consider $Y_\psi$ for different $Z$ together as different birational models of the mirror Calabi-Yau, and think of the K\"{a}hler moduli of the mirror Calabi-Yau as the union of K\"{a}hler cones of all these different birational models $Y_\psi$. This union cone is usually called the movable cone of $Y_\psi$ in algebraic geometry.\\
\begin{co}
The K\"{a}hler moduli of the mirror Calabi-Yau $Y_\psi$ can be naturally identified with\\
\[
\tau = \bigcup_{Z\in \tilde{Z}} \tau_Z
\]\\
as the set of $w=(w_m)_{m\in \Delta^0}$ (understood as a piecewise linear function on $M$) that is convex with respect to $\Delta^0$ in the sense of definition \ref{bd}, modulo linear functions on $M$.
\end{co}
\begin{flushright} $\Box$ \end{flushright}

\subsection{Monomial-divisor mirror map}
Combine the discussions from the last two subsections, we can see that the K\"{a}hler cones of various birational models of $Y_\psi$ naturally correspond to the top dimensional cones in the secondary fan for the complex moduli of $X_\psi$.\\

Mirror symmetry actually requires a precise identification of the complex moduli $P_{\tilde{\Sigma}}$ near the large complex limit and the complexified K\"{a}hler moduli $(\tilde{N}\otimes_\mathbb{Z}\mathbb{R} + i\tau)/\tilde{N}$ near the large radius limit.\\

As a toric variety, the complex moduli $P_{\tilde{\Sigma}}$ can be understood as a partial compactification of the complex torus $T_{\tilde{N}} = \tilde{N}\otimes_\mathbb{Z}\mathbb{C}/\tilde{N}$. Combine the natural embeddings $\tilde{N}\otimes_\mathbb{Z}\mathbb{R} + i\tau \subset \tilde{N}\otimes_\mathbb{Z}\mathbb{C}$ and $T_{\tilde{N}} \subset P_{\tilde{\Sigma}}$, we get the so-called {\bf monomial-divisor mirror map}

\[
\Phi: (\tilde{N}\otimes_\mathbb{Z}\mathbb{R} + i\tau)/\tilde{N} \rightarrow P_{\tilde{\Sigma}}.
\]

More precisely, for $u =(\eta_m + iw_m)_{m\in \Delta^0}\in \tilde{N}\otimes_\mathbb{Z}\mathbb{R} + i\tau$, we have the corresponding quintic polynomial\\
\[
p_u(z)= \sum_{m\in \Delta^0} e^{2\pi i\eta_m}e^{-2\pi w_m}z^m + z^{m_0}.
\]\\
(Notice $w_m$ here will correspond to $-\frac{\log t}{2\pi}w_m$ under the notation of section 2.)\\

Under the monomial-divisor mirror map, the complexified \k cone $\tau_Z^\mathbb{C} = (\tilde{N}\otimes_\mathbb{Z}\mathbb{R} + i\tau_Z)/\tilde{N}$ is natually identified with the chamber $U_Z$. In section 2, we constructed the Lagrangian torus fibration for quintic Calabi-Yau $X$ in $U_Z$ near the large complex limit. In the following sections, we will construct Lagrangian torus fibration for the mirror Calabi-Yau manifold $Y$ (more precisely the monomial-divisor map image of $X$) with \k class in $\tau_Z$ and prove the symplectic SYZ duality for the two Lagrangian torus fibrations.\\

There are several nice properties that make the secondary fan compactification a very desirable canonical partial compactification of the complex moduli of Calabi-Yau hypersurfaces in toric variety. Firstly, the secondary fan compactification is the canonical minimal partial compactification that dominates all GIT/symplectic reduction partial compactifications. Secondly, the secondary fan compactification is the minimal partial compactification such that the discriminant locus is ample. (Please consult \cite{AGM} for mentioning and references of these results.) \cite{AGM} provides the first justification that the secondary fan compactification suits the mirror symmetry purpose by identifying the top dimensional cones in the secondary fan with the \k cones of different birationally equivalent models of the mirror, therefore establishing the monomial-divisor mirror map that identifies the complexified \k moduli of the mirror with the \k cone decomposition and the complex moduli of Calabi-Yau hypersurfaces in toric variety near the large complex limit with the chamber decomposition associated with the secondary fan compactification. Our discussion in this section gives an intrinsic geometric explanation of this chamber decomposition, consequently provides direct justification for the secondary fan compactification without going to the mirror. This interpretation can potentially be generalized to determine the chamber decomposition and partial compactification of the complex moduli near the large complex limit for more general Calabi-Yau manifolds in non-toric situations. On the other hand, this interpretation also indicates the central importance of the SYZ picture in understanding mirror symmetry.\\

{\bf Remark:} The monomial-divisor mirror map first introduced in \cite{AGM} had certain restriction. Therefore, the results there can not be directly applied to our situation. As discussed in \cite{AGM}, the actual mirror map is not exactly the monomial-divisor map, or rather is a perturbation of the monomial-divisor map.\\

\se{Lagrangian torus fibrations for the mirror manifolds of quintics}
To establish the symplectic topological SYZ correspondence for quintic Calabi-Yau hypersurfaces, we need to construct Lagrangian torus fibrations for their mirror manifolds. The construction of Lagrangian torus fibrations of quintics was carried out in section 2, where the fibration structure is essentially determined by the complex structure of the quintic. In particular, the singular locus of the fibration map arises naturally from the ``string diagram'' structure of the singular set curves near the large complex limit. On the other side, the mirrors of quintics start as singular Calabi-Yau manifolds. Through crepant resolution, we get smooth models for the mirrors of quintics. Different \k structures will determine different crepant resolutions which in turn will determine different combinatoric structures of the singular locus of the Lagrangian torus fibration for the mirror of quintic Calabi-Yau manifolds. After suitable crepant resolution, the construction of Lagrangian torus fibrations for the mirror of quintics will be carried out using gradient flow. The construction is somewhat easier in this case, because most part of the singular locus is automatically 1-dimensional.\\

Take a generic $w=(w_m)_{m\in \Delta^0} \in \tau = \bigcup_{Z\in \tilde{Z}} \tau_Z$ (the K\"{a}hler moduli of mirror $Y_\psi$). In general, $w=(w_m)_{m\in \Delta^0}$ naturally defines a piecewise linear convex function $p_w$ on $M$ that satisfies $p_w(m)=-w_m$. $p_w$ naturally determines a fan $\Sigma^w$ for $M$ that is compatible with $\Delta^\vee$. For our generic choice of $w$, we may assume that the fan $\Sigma^w$ is simplicial. $p_w$ also naturally determines a real polyhedron $\Delta_w\subset N$ that consists of those $n\in N$ that as a linear function on $M$ is greater than or equal to $p_w$. As a real polyhedron, $\Delta_w$ has a dual real polyhedron $\Delta_w^\vee\subset M$.\\

\subsection{Lagrangian torus fibration with part of singular locus being codimension one}
Recall from the last section, the anti-canonical model of the mirror of quintic is the quotient of Fermat type quintic\\
\[
Y_\psi^0 = X_\psi/(\mathbb{Z}_5)^3 \subset P_{\Sigma_{\Delta^\vee}} \cong \mathbb{CP}^4/(\mathbb{Z}_5)^3.
\]\\
$w$ determines a simplical fan $\Sigma^w$ in $M$. $\hat{\pi}: P_{\Sigma^w} \rightarrow P_{\Sigma_{\Delta^\vee}}$ is a resolution. The pullback $Y_\psi = \hat{\pi}^{-1}(Y_\psi^0)$ is the mirror Calabi-Yau with K\"{a}hler class corresponding to $w=(w_m)_{m\in \Delta^0}$.\\

On $P_{\Sigma^w}$ we may take any toric metric that represents the \k class of $w$. For example, we may use the toric metric resulting from symplectic reduction as discussed in Guillemin's paper \cite{Gu}, whose \k potential $h_w$ is the Legendre transformation of

\[
v_w(n) = \sum_{m\in \Delta^0} (\langle m,n\rangle +w_m)(\log (\langle m,n\rangle +w_m) -1).
\]

Let $\Delta_w^{(k)}$ denote the smooth open part of the $k$-skeleton of $\Delta_w$. Similarly, let $Y_\infty^{(k)}$ denote the smooth open part of the complex $k$-skeleton of $Y_\infty$. Then we have the stratification

\[
\Delta_w = \bigcup_{k=0}^3 \Delta_w^{(k)},\ \ \ Y_\infty = \bigcup_{k=0}^3 Y_\infty^{(k)}.
\]

A torus fibration $F_\infty: Y_\infty \rightarrow \partial \Delta_w$ is called topologically smooth if $\Delta_w^{(k)} = F_\infty(Y_\infty^{(k)})$ and $F_\infty|_{Y_\infty^{(k)}}: Y_\infty^{(k)} \rightarrow \Delta_w^{(k)}$ is a topologically smooth fibration, whose fibers are necessarily $k$-dimensional torus.\\

\begin{prop}
\label{ff}
The moment map $F_w$ maps $P_{\Sigma^w}$ to $\Delta_w$, and maps $Y_\infty \subset P_{\Sigma^w}$ to $\partial \Delta_w$. $F_w|_{Y_\infty}: Y_\infty \rightarrow \partial \Delta_w$ is a topologically smooth Lagrangian torus fibration with respect to the symplectic form $\omega = \partial \bar{\partial}\log h_w$ on $P_{\Sigma^w}$.
\end{prop}
\begin{flushright} $\Box$ \end{flushright}
With the Lagrangian fibration of $Y_\infty$ in hand, just as in the case of quintics, we can similarly apply theorem 8.1 in \cite{lag2} to produce the Lagrangian torus fibration $F_\psi: Y_\psi \rightarrow \partial \Delta_w$ for $Y_\psi$. As in the quintic case, we will rephrase the general theorem 8.1 of \cite{lag2} in the situation of mirror of quintic.\\

\begin{theorem}
\label{fd}
Start with a topologically smooth Lagrangian torus fibration $F_\infty: Y_\infty \rightarrow \partial \Delta_w$, we can construct a symplectic morphism $H_\psi: Y_\psi \rightarrow Y_\infty$ such that $F_\psi = F_\infty\circ H_\psi: Y_\psi \rightarrow \partial \Delta_w$ is a Lagrangian torus fibration with singular set $C = Y_\psi \cap {\rm Sing}(Y_\infty)$ and singular locus $\Gamma = F_\infty(C)$. For $b\not\in \Gamma$, $F_\psi^{-1}(b)$ is a real $3$-torus. For $b \in \Gamma$, $F_\psi^{-1}(b)$ is singular. For $b \in \Gamma \cap \Delta_w^{(k)}$, $F_\psi^{-1}(b)\cap C = F_\infty^{-1}(b)\cap C$ and $H_\psi: F_\psi^{-1}(b)\backslash C \rightarrow F_\infty^{-1}(b)\backslash C$ is a topologically smooth $(3-k)$-torus fibration.
\end{theorem}
\begin{flushright} $\Box$ \end{flushright}
We will start with the natural $F_\infty = F_w|_{Y_\infty}$. The task now is to understand the structure of singular fibers and the singular locus of the corresponding fibration $F_\psi$. The singular locus $\tilde{\Gamma} = F_\psi(C)$ is a fattening of some graph $\Gamma$. To describe this graph $\Gamma$, let us recall from the last section of \cite{lag1} that the singular locus of the Lagrangian fibration of $Y^0_\psi \subset P_{\Sigma_{\Delta^\vee}}$ is a fattening of a graph $\hat{\Gamma}\subset \partial \Delta^\vee$, where\\
\[
\hat{\Gamma} = \bigcup_{\{ijklm\}=\{12345\}}\overline{P_{ij}P_{klm}}.
\]\\
The following is a picture (from \cite{lag1}) of a face of $\partial \Delta^\vee$.
\begin{center}
\leavevmode
\hbox{%
\epsfxsize=4in
\epsffile{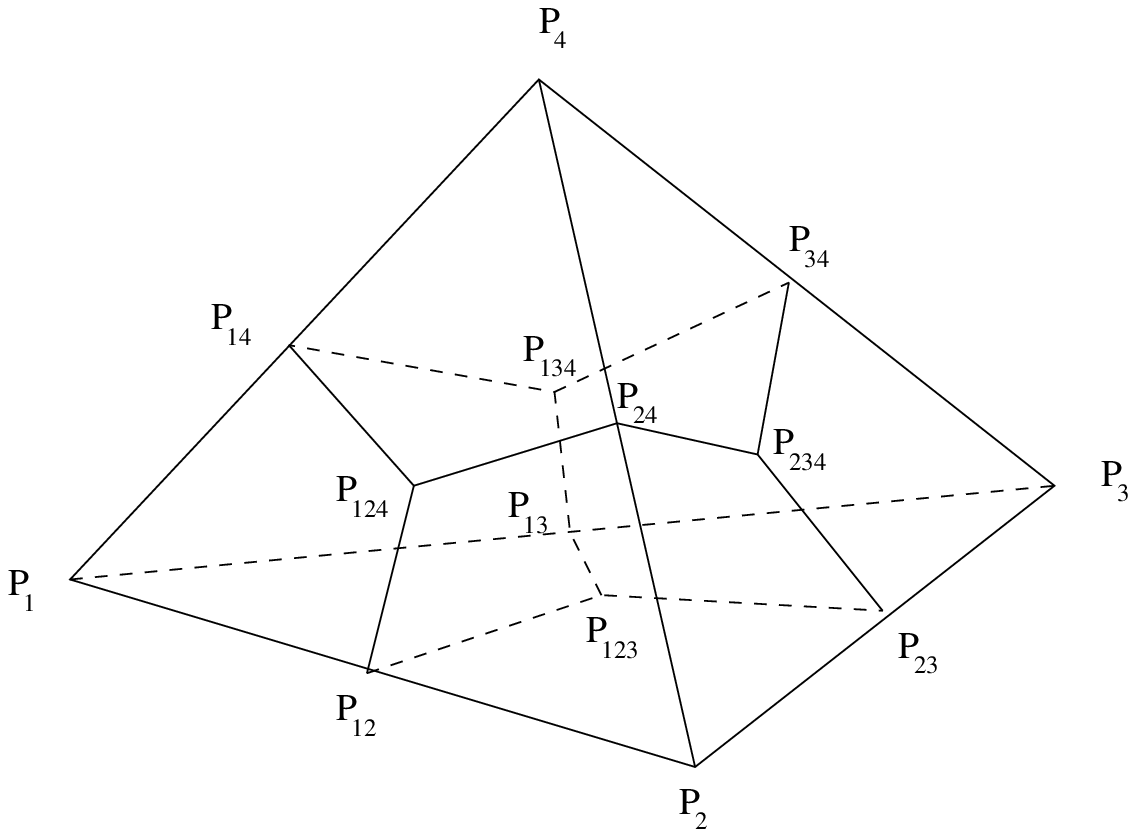}}
\end{center}
\begin{center}
\stepcounter{figure}
Figure \thefigure: $\hat{\Gamma}$ in a face of $\partial \Delta^\vee$
\end{center}
It is interesting to observe that\\
\[
{\rm Sing}(P_{\Sigma_{\Delta^\vee}}) = {\rm Sing}(Y_\infty).
\]
Hence\\
\[
{\rm Sing}(Y^0_\psi) = C = Y^0_\psi \cap {\rm Sing}(Y_\infty).
\]
Let $\tilde{P}_{ij}$ be the point in $C$ that maps to $P_{ij}$ in $\Delta^\vee$. Notice that ${\rm Sing}(C) = \{\tilde{P}_{ij}\}$. Along the smooth part of $C$, $Y^0_\psi$ has $A_5$-singularity. Under the unique crepant resolution, $C$ is turned into 5 copies of $C$. Around $\tilde{P}_{ij} \in {\rm Sing}(C)$, singularity of $Y^0_\psi$ is much more complicated and crepant resolution is not unique. The following is a picture (from \cite{lag1}) of the fan of such singularity and the subdivision fan of the standard crepant resolution.\\
\begin{center}
\leavevmode
\hbox{%
\epsfxsize=5in
\epsffile{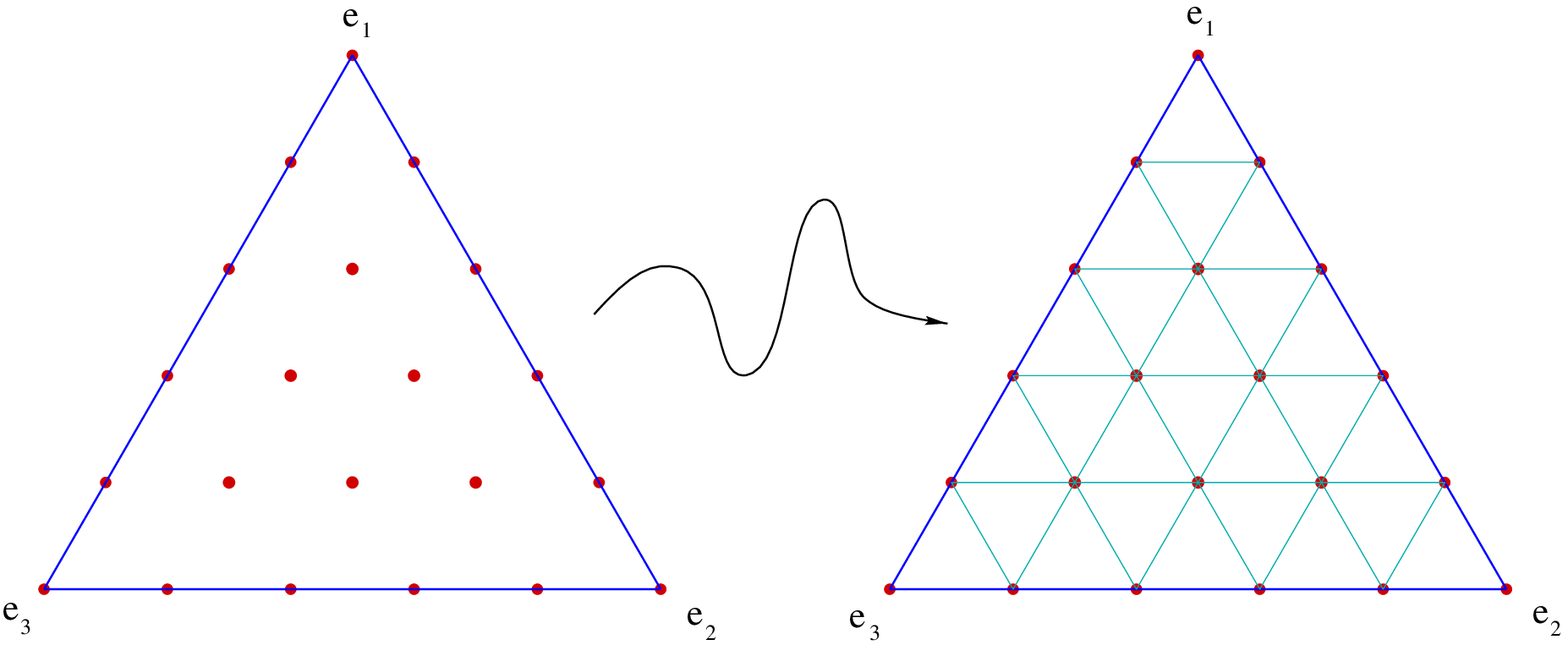}}
\end{center}
\begin{center}
\stepcounter{figure}
Figure \thefigure: the standard crepant resolution of singularity at $\hat{P}_{ij}$
\end{center}
$\hat{\pi}: P_{\Sigma^w} \rightarrow P_{\Sigma_{\Delta^\vee}}$ naturally induces a map $\pi: \Delta_w \rightarrow \Delta^\vee$. We may take $\Gamma$ to be the 1-skeleton of $\pi^{-1}(\hat{\Gamma})$. In a small neighborhood of $P_{ij}$, $\hat{\Gamma} \subset \partial \Delta^\vee$ is indicated in the following picture (Figure 7):\\
\begin{center}
\setlength{\unitlength}{1pt}
\begin{picture}(200,100)(-30,10)
\put(80,80){\circle*{2.5}}
\put(90,82){\circle*{2}}
\put(70,78){\circle*{2}}

\multiput(80,80)(36,0){1}{\line(5,1){10}}
\multiput(80,17)(23,0){1}{\line(5,1){10}}
\multiput(26,107)(72,0){1}{\line(5,1){10}}
\multiput(134,107)(36,0){1}{\line(5,1){10}}
\multiput(80,80)(36,0){1}{\line(-5,-1){10}}
\multiput(80,17)(23,0){1}{\line(-5,-1){10}}
\multiput(26,107)(72,0){1}{\line(-5,-1){10}}
\multiput(134,107)(36,0){1}{\line(-5,-1){10}}
\multiput(90,82)(36,0){1}{\line(2,1){54}}
\multiput(90,82)(36,0){1}{\line(-2,1){54}}
\multiput(90,82)(36,0){1}{\line(0,-1){63}}
\multiput(70,78)(36,0){1}{\line(2,1){54}}
\multiput(70,78)(36,0){1}{\line(-2,1){54}}
\multiput(70,78)(36,0){1}{\line(0,-1){63}}
\put(102,48){\vector(-1,0){21}}
\put(110,50){\makebox(0,0){$\hat{\Gamma}$}}
\thicklines
\multiput(80,80)(36,0){1}{\line(2,1){54}}
\multiput(80,80)(36,0){1}{\line(-2,1){54}}
\multiput(80,80)(36,0){1}{\line(0,-1){63}}
\end{picture}
\end{center}
\begin{center}
\stepcounter{figure}
Figure \thefigure: $\hat{\Gamma}\subset \partial \Delta^\vee$ near $P_{ij}$
\end{center}
Under the standard crepant resolution, we get $\Gamma \subset \partial \Delta_w$ as indicated in the following picture (Figure 8):\\
\begin{center}
\setlength{\unitlength}{1pt}
\begin{picture}(200,160)(-30,0)
\multiput(78,149)(36,0){1}{\line(5,1){10}}
\multiput(60,119)(36,0){2}{\line(5,1){10}}
\multiput(42,89)(36,0){3}{\line(5,1){10}}
\multiput(24,59)(36,0){4}{\line(5,1){10}}
\multiput(6,29)(36,0){5}{\line(5,1){10}}

\multiput(78,128)(36,0){1}{\line(5,1){10}}
\multiput(60,98)(36,0){2}{\line(5,1){10}}
\multiput(42,68)(36,0){3}{\line(5,1){10}}
\multiput(24,38)(36,0){4}{\line(5,1){10}}

\multiput(6,8)(36,0){5}{\line(5,1){10}}
\multiput(-12,38)(18,30){5}{\line(5,1){10}}
\multiput(96,158)(18,-30){5}{\line(5,1){10}}

\multiput(78,149)(36,0){1}{\line(-5,-1){10}}
\multiput(60,119)(36,0){2}{\line(-5,-1){10}}
\multiput(42,89)(36,0){3}{\line(-5,-1){10}}
\multiput(24,59)(36,0){4}{\line(-5,-1){10}}
\multiput(6,29)(36,0){5}{\line(-5,-1){10}}

\multiput(78,128)(36,0){1}{\line(-5,-1){10}}
\multiput(60,98)(36,0){2}{\line(-5,-1){10}}
\multiput(42,68)(36,0){3}{\line(-5,-1){10}}
\multiput(24,38)(36,0){4}{\line(-5,-1){10}}

\multiput(6,8)(36,0){5}{\line(-5,-1){10}}
\multiput(-12,38)(18,30){5}{\line(-5,-1){10}}
\multiput(96,158)(18,-30){5}{\line(-5,-1){10}}
\thicklines
\multiput(78,149)(36,0){1}{\line(2,1){18}}
\multiput(78,149)(36,0){1}{\line(-2,1){18}}
\multiput(78,128)(36,0){1}{\line(0,1){21}}
\multiput(42,128)(36,0){2}{\line(2,-1){18}}
\multiput(78,128)(36,0){2}{\line(-2,-1){18}}
\multiput(60,98)(36,0){2}{\line(0,1){21}}
\multiput(24,98)(36,0){3}{\line(2,-1){18}}
\multiput(60,98)(36,0){3}{\line(-2,-1){18}}
\multiput(42,68)(36,0){3}{\line(0,1){21}}
\multiput(6,68)(36,0){4}{\line(2,-1){18}}
\multiput(42,68)(36,0){4}{\line(-2,-1){18}}
\multiput(24,38)(36,0){4}{\line(0,1){21}}
\multiput(-12,38)(36,0){5}{\line(2,-1){18}}
\multiput(24,38)(36,0){5}{\line(-2,-1){18}}
\multiput(6,8)(36,0){5}{\line(0,1){21}}

\thinlines
\put(10,2){
\multiput(78,149)(36,0){1}{\line(2,1){18}}
\multiput(78,149)(36,0){1}{\line(-2,1){18}}
\multiput(78,128)(36,0){1}{\line(0,1){21}}
\multiput(42,128)(36,0){2}{\line(2,-1){18}}
\multiput(78,128)(36,0){2}{\line(-2,-1){18}}
\multiput(60,98)(36,0){2}{\line(0,1){21}}
\multiput(24,98)(36,0){3}{\line(2,-1){18}}
\multiput(60,98)(36,0){3}{\line(-2,-1){18}}
\multiput(42,68)(36,0){3}{\line(0,1){21}}
\multiput(6,68)(36,0){4}{\line(2,-1){18}}
\multiput(42,68)(36,0){4}{\line(-2,-1){18}}
\multiput(24,38)(36,0){4}{\line(0,1){21}}
\multiput(-12,38)(36,0){5}{\line(2,-1){18}}
\multiput(24,38)(36,0){5}{\line(-2,-1){18}}
\multiput(6,8)(36,0){5}{\line(0,1){21}}}

\put(-10,-2){
\multiput(78,149)(36,0){1}{\line(2,1){18}}
\multiput(78,149)(36,0){1}{\line(-2,1){18}}
\multiput(78,128)(36,0){1}{\line(0,1){21}}
\multiput(42,128)(36,0){2}{\line(2,-1){18}}
\multiput(78,128)(36,0){2}{\line(-2,-1){18}}
\multiput(60,98)(36,0){2}{\line(0,1){21}}
\multiput(24,98)(36,0){3}{\line(2,-1){18}}
\multiput(60,98)(36,0){3}{\line(-2,-1){18}}
\multiput(42,68)(36,0){3}{\line(0,1){21}}
\multiput(6,68)(36,0){4}{\line(2,-1){18}}
\multiput(42,68)(36,0){4}{\line(-2,-1){18}}
\multiput(24,38)(36,0){4}{\line(0,1){21}}
\multiput(-12,38)(36,0){5}{\line(2,-1){18}}
\multiput(24,38)(36,0){5}{\line(-2,-1){18}}
\multiput(6,8)(36,0){5}{\line(0,1){21}}}
\put(172,18){\vector(-1,0){21}}
\put(179,18){\makebox(0,0){$\Gamma$}}
\end{picture}\\
\end{center}
\begin{center}
\stepcounter{figure}
Figure \thefigure: $\Gamma$ for the standard crepant resolution
\end{center}
For a different crepant resolution, we can get an alternative picture for $\Gamma \subset \partial \Delta_w$ (Figure 9).\\
\begin{center}
\setlength{\unitlength}{1pt}
\begin{picture}(200,160)(-30,0)
\multiput(78,149)(36,0){1}{\line(5,1){10}}
\multiput(60,119)(23,0){2}{\line(5,1){10}}
\multiput(42,89)(72,0){2}{\line(5,1){10}}
\multiput(24,59)(36,0){4}{\line(5,1){10}}
\multiput(6,29)(36,0){5}{\line(5,1){10}}

\multiput(91,128)(36,0){1}{\line(5,1){10}}
\multiput(60,98)(36,0){2}{\line(5,1){10}}
\multiput(42,68)(72,0){2}{\line(5,1){10}}
\multiput(24,38)(36,0){4}{\line(5,1){10}}
\multiput(69.75,78.5)(16.5,0){2}{\line(5,1){10}}

\multiput(6,8)(36,0){5}{\line(5,1){10}}
\multiput(-12,38)(18,30){5}{\line(5,1){10}}
\multiput(96,158)(18,-30){5}{\line(5,1){10}}

\multiput(78,149)(36,0){1}{\line(-5,-1){10}}
\multiput(60,119)(23,0){2}{\line(-5,-1){10}}
\multiput(42,89)(72,0){2}{\line(-5,-1){10}}
\multiput(24,59)(36,0){4}{\line(-5,-1){10}}
\multiput(6,29)(36,0){5}{\line(-5,-1){10}}

\multiput(91,128)(36,0){1}{\line(-5,-1){10}}
\multiput(60,98)(36,0){2}{\line(-5,-1){10}}
\multiput(42,68)(72,0){2}{\line(-5,-1){10}}
\multiput(24,38)(36,0){4}{\line(-5,-1){10}}
\multiput(69.75,78.5)(16.5,0){2}{\line(-5,-1){10}}

\multiput(6,8)(36,0){5}{\line(-5,-1){10}}
\multiput(-12,38)(18,30){5}{\line(-5,-1){10}}
\multiput(96,158)(18,-30){5}{\line(-5,-1){10}}
\thicklines
\multiput(78,149)(36,0){1}{\line(2,1){18}}
\multiput(78,149)(36,0){1}{\line(-2,1){18}}
\multiput(78,149)(4,-30){2}{\line(2,-3){14}}
\multiput(114,128)(-31,-9){2}{\line(-1,0){23}}
\multiput(82,119)(-27,-10){1}{\line(1,1){10}}
\multiput(42,128)(36,0){1}{\line(2,-1){18}}
\multiput(78,128)(36,0){0}{\line(-2,-1){18}}
\multiput(60,98)(36,0){1}{\line(0,1){21}}
\multiput(24,98)(72,0){2}{\line(2,-1){18}}
\multiput(60,98)(26.5,-19.5){2}{\line(1,-2){9.75}}
\put(69.75,78.5){\line(1,0){16.5}}
\multiput(60,98)(72,0){2}{\line(-2,-1){18}}
\multiput(96,98)(-26.5,-19.5){2}{\line(-1,-2){9.75}}
\multiput(42,68)(72,0){2}{\line(0,1){21}}
\multiput(6,68)(36,0){2}{\line(2,-1){18}}
\multiput(114,68)(36,0){1}{\line(2,-1){18}}
\multiput(42,68)(36,0){1}{\line(-2,-1){18}}
\multiput(114,68)(36,0){2}{\line(-2,-1){18}}
\multiput(24,38)(36,0){4}{\line(0,1){21}}
\multiput(-12,38)(36,0){5}{\line(2,-1){18}}
\multiput(24,38)(36,0){5}{\line(-2,-1){18}}
\multiput(6,8)(36,0){5}{\line(0,1){21}}

\thinlines
\put(10,2){
\multiput(78,149)(36,0){1}{\line(2,1){18}}
\multiput(78,149)(36,0){1}{\line(-2,1){18}}
\multiput(78,149)(4,-30){2}{\line(2,-3){14}}
\multiput(114,128)(-31,-9){2}{\line(-1,0){23}}
\multiput(82,119)(-27,-10){1}{\line(1,1){10}}
\multiput(42,128)(36,0){1}{\line(2,-1){18}}
\multiput(78,128)(36,0){0}{\line(-2,-1){18}}
\multiput(60,98)(36,0){1}{\line(0,1){21}}
\multiput(24,98)(72,0){2}{\line(2,-1){18}}
\multiput(60,98)(26.5,-19.5){2}{\line(1,-2){9.75}}
\put(69.75,78.5){\line(1,0){16.5}}
\multiput(60,98)(72,0){2}{\line(-2,-1){18}}
\multiput(96,98)(-26.5,-19.5){2}{\line(-1,-2){9.75}}
\multiput(42,68)(72,0){2}{\line(0,1){21}}
\multiput(6,68)(36,0){2}{\line(2,-1){18}}
\multiput(114,68)(36,0){1}{\line(2,-1){18}}
\multiput(42,68)(36,0){1}{\line(-2,-1){18}}
\multiput(114,68)(36,0){2}{\line(-2,-1){18}}
\multiput(24,38)(36,0){4}{\line(0,1){21}}
\multiput(-12,38)(36,0){5}{\line(2,-1){18}}
\multiput(24,38)(36,0){5}{\line(-2,-1){18}}
\multiput(6,8)(36,0){5}{\line(0,1){21}}}

\put(-10,-2){
\multiput(78,149)(36,0){1}{\line(2,1){18}}
\multiput(78,149)(36,0){1}{\line(-2,1){18}}
\multiput(78,149)(4,-30){2}{\line(2,-3){14}}
\multiput(114,128)(-31,-9){2}{\line(-1,0){23}}
\multiput(82,119)(-27,-10){1}{\line(1,1){10}}
\multiput(42,128)(36,0){1}{\line(2,-1){18}}
\multiput(78,128)(36,0){0}{\line(-2,-1){18}}
\multiput(60,98)(36,0){1}{\line(0,1){21}}
\multiput(24,98)(72,0){2}{\line(2,-1){18}}
\multiput(60,98)(26.5,-19.5){2}{\line(1,-2){9.75}}
\put(69.75,78.5){\line(1,0){16.5}}
\multiput(60,98)(72,0){2}{\line(-2,-1){18}}
\multiput(96,98)(-26.5,-19.5){2}{\line(-1,-2){9.75}}
\multiput(42,68)(72,0){2}{\line(0,1){21}}
\multiput(6,68)(36,0){2}{\line(2,-1){18}}
\multiput(114,68)(36,0){1}{\line(2,-1){18}}
\multiput(42,68)(36,0){1}{\line(-2,-1){18}}
\multiput(114,68)(36,0){2}{\line(-2,-1){18}}
\multiput(24,38)(36,0){4}{\line(0,1){21}}
\multiput(-12,38)(36,0){5}{\line(2,-1){18}}
\multiput(24,38)(36,0){5}{\line(-2,-1){18}}
\multiput(6,8)(36,0){5}{\line(0,1){21}}}
\put(172,18){\vector(-1,0){21}}
\put(179,18){\makebox(0,0){$\Gamma$}}
\end{picture}
\end{center}
\begin{center}
\stepcounter{figure}
Figure \thefigure: $\Gamma$ for alternative crepant resolution
\end{center}
These singular locus graphs clearly resemble the singular locus graphs in section 2 via string diagram (amoeba) construction, although the two constructions are quite different.\\
\begin{prop}
\label{fc}
There exist a piecewise smooth automorphism of $\Delta_w$ that preserves subfaces of $\Delta_w$ so that under such automorphism, the 1-skeleton $\Gamma$ of $\pi^{-1}(\hat{\Gamma})$ is identified with the union of 1-simplices in the baricenter subdivision of the 2-skeleton of $\Delta_w$ that under the map $\pi$ do not meet the vertices of $\Delta^\vee$.
\end{prop}
{\bf Proof:} Clearly, one only need to construct such automorphism for each of the 2-subfaces of $\Delta_w$ that contain part of $\Gamma$ and then piecewise smoothly extend to the whole $\Delta_w$. There are only two kinds of such 2-subfaces, those map to 1-subsimpleces or 2-subsimpleces in $\Delta^\vee$ under $\pi$. In both cases the required automorphisms are rather trivial to construct.
\begin{flushright} $\Box$ \end{flushright}
Let $\tilde{\Gamma}_0$ denote the singular locus for the fibration of $Y^0_\psi \subset P_{\Sigma_{\Delta^\vee}}$. Similar to the decomposition (\ref{bg}) in the quintic case, we have

\[
\tilde{\Gamma}_0 = \tilde{\Gamma}^0_0 \cup \tilde{\Gamma}^1_0 \cup \tilde{\Gamma}^2_0.
\]

Consider the map $\pi|_{\tilde{\Gamma}}: \tilde{\Gamma} \rightarrow \tilde{\Gamma}_0$. For $k=1,2$, define $\tilde{\Gamma}^k = (\pi|_{\tilde{\Gamma}})^{-1}(\tilde{\Gamma}^k_0)$. Let $\Gamma' = (\pi|_{\tilde{\Gamma}})^{-1}(\tilde{\Gamma}^0_0)$. $\Gamma'$ represents the part of $\Gamma$ that is not fattened in $\tilde{\Gamma}$. $\Gamma'$ has the natural decomposition $\Gamma' = (\Gamma')^1\cup \Gamma^3$, where $(\Gamma')^1$ is the smooth part of $\Gamma'$, $\Gamma^3$ is the singular part of $\Gamma$ in the 1-skeleton of $\Delta_w$. Together we have the decomposition

\[
\tilde{\Gamma} = (\Gamma')^1\cup \Gamma^3 \cup \tilde{\Gamma}^1 \cup \tilde{\Gamma}^2.
\]

Start with the natural $F_\infty = F_w|_{Y_\infty}$ in proposition \ref{ff}, apply theorem \ref{fd}, we have

\begin{theorem}
For $Y_\psi \subset P_{\Sigma^w}$, when $w=(w_m)_{m\in \Delta^0}\in \tau$ is generic and near the large radius limit of $\tau$, the gradient flow method will produce a Lagrangian fibration $F_\psi: Y_{\psi} \rightarrow \partial\Delta_w$ with singular locus $\tilde{\Gamma} = F_\infty(C)$ as fattening of graph $\Gamma$. There are five types of fibers.\\
(i). For $p\in \partial\Delta_w\backslash \tilde{\Gamma}$, $F_\psi^{-1}(p)$ is a smooth Lagrangian 3-torus.\\
(ii). For $p\in \tilde{\Gamma}^2$, $F_\psi^{-1}(p)$ is a Lagrangian 3-torus with two circles collapsed to two singular points.\\
(iii). For $p\in \tilde{\Gamma}^1$, $F_\psi^{-1}(p)$ is a Lagrangian 3-torus with one circle collapsed to one singular point.\\
(iv). For $p\in \Gamma^3$, $F_{\psi}^{-1}(p)$ is a type III singular fiber.\\
(v). For $p\in (\Gamma')^1$, $F_{\psi}^{-1}(p)$ is a type I singular fiber.
\end{theorem}
\begin{flushright} $\Box$ \end{flushright}
{\bf Remark:} From the construction of this Lagrangian fibration $F_\psi$, one can already see the hint of SYZ mirror correspondence to the Lagrangian fibrations of quintics. Although the singular locus $\tilde{\Gamma}$ of $F_\psi$ is of codimension one, $\tilde{\Gamma}$ is only fattening the part of $\Gamma$ which is mapped to the interior of 2-simplices of $\Delta^\vee$ under the map $\pi$. A significant part $\Gamma'$ of $\Gamma$ that results from the crepant resolution is not fattened. The SYZ mirror correspondence is already quite apparent on this part of the fibration.\\

\subsection{Non-Lagrangian torus fibration with codimension two singular locus}
The commuting diagram
\begin{center}
\setlength{\unitlength}{1.2pt}
\begin{picture}(100,50)(0,0)
\put(30,40){\makebox(0,0){$P_{\Sigma^w}$}}
\put(80,40){\makebox(0,0){$P_{\Sigma_{\Delta^\vee}}$}}
\put(30,10){\makebox(0,0){$\Delta_w$}}
\put(80,10){\makebox(0,0){$\Delta^\vee$}}
\put(42,40){\vector(1,0){26}}
\put(28,32){\vector(0,-1){15}}
\put(42,10){\vector(1,0){26}}
\put(80,32){\vector(0,-1){15}}
\put(55,45){\makebox(0,0){\footnotesize{$\hat{\pi}$}}}
\put(55,15){\makebox(0,0){\footnotesize{$\pi$}}}
\put(38,25){\makebox(0,0){\footnotesize{$F_w$}}}
\put(85,25){\makebox(0,0){\footnotesize{$F$}}}
\end{picture}
\end{center}
gives rise to
\begin{center}
\setlength{\unitlength}{1.2pt}
\begin{picture}(100,50)(0,0)
\put(30,40){\makebox(0,0){$Y_\infty$}}
\put(80,40){\makebox(0,0){$Y^0_\infty$}}
\put(30,10){\makebox(0,0){$\partial\Delta_w$}}
\put(80,10){\makebox(0,0){$\partial\Delta^\vee$}}
\put(42,40){\vector(1,0){26}}
\put(28,32){\vector(0,-1){15}}
\put(42,10){\vector(1,0){26}}
\put(80,32){\vector(0,-1){15}}
\put(55,45){\makebox(0,0){\footnotesize{$\hat{\pi}$}}}
\put(55,15){\makebox(0,0){\footnotesize{$\pi$}}}
\put(38,25){\makebox(0,0){\footnotesize{$F_w$}}}
\put(85,25){\makebox(0,0){\footnotesize{$F$}}}
\end{picture}
\end{center}
Let $C^0 = Y^0_\psi \cap {\rm Sing}(Y^0_\infty)$. In \cite{lag1}, the Lagrangian fibration of $X_\infty$ is modified to yield Lagrangian torus fibration with codimension two singular locus for $X_\psi$. Notice that the modified Lagrangian fibration of $X_\infty$ and the gradient flow of the Fermat type quintic Calabi-Yau family $\{X_\psi\}$ are invariant under the action of $(\mathbb{Z}_5)^3$. The quotient gives us the corresponding gradient flow on $P_{\Sigma_{\Delta^\vee}} \cong \mathbb{CP}^4/(\mathbb{Z}_5)^3$ of the family $\{Y^0_\psi\}$, which produces the Lagrangian torus fibration with codimension 2 singular locus for $Y^0_\psi$ under the quotient orbifold \k metric. Here we will only make use of the modified torus fibration $\hat{F}: Y^0_\infty \rightarrow \partial\Delta^\vee$ that satisfies $\hat{F}(C^0) = \hat{\Gamma}$.\\
\begin{lm}
There exists torus fibration $\hat{F}_w: Y_\infty \rightarrow \partial\Delta_w$ that makes the following diagram commute.
\begin{center}
\setlength{\unitlength}{1.2pt}
\begin{picture}(100,50)(0,0)
\put(30,40){\makebox(0,0){$Y_\infty$}}
\put(80,40){\makebox(0,0){$Y^0_\infty$}}
\put(30,10){\makebox(0,0){$\partial\Delta_w$}}
\put(80,10){\makebox(0,0){$\partial\Delta^\vee$}}
\put(42,40){\vector(1,0){26}}
\put(28,32){\vector(0,-1){15}}
\put(42,10){\vector(1,0){26}}
\put(80,32){\vector(0,-1){15}}
\put(55,45){\makebox(0,0){\footnotesize{$\hat{\pi}$}}}
\put(55,15){\makebox(0,0){\footnotesize{$\pi$}}}
\put(38,25){\makebox(0,0){\footnotesize{$\hat{F}_w$}}}
\put(85,25){\makebox(0,0){\footnotesize{$\hat{F}$}}}
\end{picture}
\end{center}
\end{lm}
{\bf Proof:} For any $x\in \partial\Delta^\vee$ and $y\in \hat{F}^{-1}(x)$, naturally $F_w: \hat{\pi}^{-1}(y) \rightarrow \pi^{-1}(F(y))$. Observe that $\hat{F}$ respects the stratification of $Y_\infty$ and $\partial\Delta_w$. $F(y)$ and $x$ are in the same strata. There is a natural identification between $\pi^{-1}(F(y))$ and $\pi^{-1}(x)$. Compose this identification with $F_w$, we get $\hat{F}_w: \hat{\pi}^{-1}(y) \rightarrow \pi^{-1}(x)$. Piece together the definition for all $(x,y)$, we get the desired $\hat{F}_w$.
\begin{flushright} $\Box$ \end{flushright}
{\bf Remark:} $\hat{F}_w$ so constructed is not Lagrangian fibration in general.\\

For $C = Y_\psi \cap {\rm Sing}(Y_\infty)$, the $\hat{F}_w$ so constructed satisfies that $\hat{F}_w(C) = \Gamma$ is a graph. Let $\Gamma = \Gamma^1\cup \Gamma^2\cup \Gamma^3$, where $\Gamma^1$ is the smooth part of $\Gamma$, $\Gamma^3$ is the singular part of $\Gamma$ in the 1-skeleton of $\Delta_w$, $\Gamma^2$ is the rest of singular part of $\Gamma$. Apply theorem \ref{fd} without the Lagrangian condition, we have\\
\begin{theorem}
For $Y_\psi \subset P_{\Sigma^w}$, when $w=(w_m)_{m\in \Delta^0}\in \tau$ is generic and near the large radius limit of $\tau$, start with torus fibration $\hat{F}_w$ the gradient flow method will produce a torus fibration $\hat{F}_{\psi}: Y_{\psi} \rightarrow \partial\Delta_w$ with codimension 2 singular locus $\Gamma = \hat{F}_w(C)$. There are 4 types of fibers.\\
(i). For $p\in \partial\Delta_w\backslash \Gamma$, $\hat{F}_{\psi}^{-1}(p)$ is a 3-torus.\\
(ii). For $p\in \Gamma^1$, $\hat{F}_{\psi}^{-1}(p)$ is a type I singular fiber.\\
(iii). For $p\in \Gamma^2$, $\hat{F}_{\psi}^{-1}(p)$ is a type II singular fiber.\\
(iv). For $p\in \Gamma^3$, $\hat{F}_{\psi}^{-1}(p)$ is a type III singular fiber.
\end{theorem}
\begin{flushright} $\Box$ \end{flushright}
{\bf Remark:} Lagrangian torus fibration is much more difficult to construct than non-Lagrangian torus fibration. If our purpose is merely to construct non-Lagrangian fibration, the kind of non-Lagrangian torus fibration in the above theorem can be constructed with much simpler method. Of course to get an elegant and canonical construction is another matter.\\

\subsection{Lagrangian torus fibration with codimension two singular locus}
Take a generic $w=(w_m)_{m\in \Delta^0} \in \tau = \bigcup_{Z\in \tilde{Z}} \tau_Z$ (the K\"{a}hler moduli of mirror $Y_\psi$). Since $\Delta$ is a simplex, without loss of generality, we can normalize $w$ to assume that $w_m = w^0$ are all the same for $m$ being the vertices of $\Delta$. It is easy to see that $\Delta_w\subset w^0\Delta^\vee\subset N$, and $\Delta_w$ geometrically can be viewed as $w^0\Delta^\vee$ with some corners chopped off.\\

According to the proof of theorem \ref{be}, the construction of Lagrangian torus fibration with codimension two singular locus can be reduced to establishing a symplectic isotopy from $C = Y_\psi \cap {\rm Sing}(Y_\infty)$ to $\hat{C}$ such that $F_w(\hat{C}) = \Gamma$ is a graph. As we know, $\tilde{\Gamma} = F_w(C)$ is mostly graph. The non-graph part of $\tilde{\Gamma}$ is a disjoint union of many curved triangle pieces (50 pieces to be exact). Each of such curved triangle is mapped via $\pi$ into a 2-subsimplex in $\Delta^\vee$. Each of such 2-subsimplex corresponds to a sub-$\mathbb{CP}^2$ in $P_{\Sigma_{\Delta^\vee}}$.\\

Let us concentrate on one of these $\mathbb{CP}^2$'s. For convenience of notation, we will temporarily suspend all our notation convention and use $\Delta$ and $\Delta_w$ to denote the moment map ($F$ and $F_w$) images of $(\mathbb{CP}^2,\omega_{FS})$ and its toric blow up $(\mathbb{C}\hat{\mathbb{P}}^2,\omega_w)$. ($\Delta_w$ geometrically can be viewed as $w^0\Delta$ chopping off some corners.) We will use $C_0$ to denote both the curve $\{z_1+z_2+1=0\} \subset (\mathbb{CP}^2,\omega_{FS})$ and its proper transformation in $(\mathbb{C}\hat{\mathbb{P}}^2,\omega_w)$. In such way, our construction is reduced to the following problem on $\mathbb{CP}^2$. \\

{\bf Problem:} Find a symplectic isotopy from $C_0$ to symplectic curve $\hat{C}\subset (\mathbb{C}\hat{\mathbb{P}}^2,\omega_w)$ such that $F_w(\hat{C})$ is a ``Y'' shaped graph in $\Delta_w$.\\

Naturally $\Delta_w \subset w^0\Delta$. Introduce the family of toric \k forms $\omega_t = (1-t)w^0\omega_{FS} + t\omega_w$. Then we have the family $F_t: (\mathbb{C}\hat{\mathbb{P}}^2,\omega_t) \rightarrow \Delta_t$. $F_t =(1-t)w^0F + tF_w$.  Assume $F_0(C_0)\subset \Delta_w$. Since $F_w(C_0)\subset \Delta_w$, we have $F_t(C_0)\subset \Delta_w$ for all $0\leq t \leq 1$. Take a neighborhood $U$ of $\cup_t F_t(C_0)$ in $\Delta_w$. Then we have:

\begin{lm}
The natural maps $f_t: (F_t^{-1}(U),\omega_t) \rightarrow (F_1^{-1}(U),\omega_1)$ are symplectomorphisms.
\end{lm}
{\bf Proof:} Assume $F_t = (h_1^t, h_2^t)$, then

\[
\omega_t = d\theta_1\wedge dh_1^t + d\theta_2\wedge dh_2^t.
\]

The natural map $f_t$ satisfies $F_1\circ f_t = F_t$. Namely $h_i^1\circ f_t = h_i^t$. Therefore $f_t^*\omega_1 = \omega_t$.
\begin{flushright} $\Box$ \end{flushright}
\begin{prop}
There exists a symplectic isotopy from $C_0$ to symplectic curve $\hat{C}\subset (\mathbb{C}\hat{\mathbb{P}}^2,\omega_w)$ such that $F_w(\hat{C})$ is a ``Y'' shaped graph in $\Delta_w$.
\end{prop}
{\bf Proof:} According to theorem 4.1 in \cite{N}, there exists a symplectic isotopy $\{C_t\}_{t\in[0,1]}$ from $C_0$ to $\hat{C}_0$ in $\mathbb{CP}^2$ such that $F(\hat{C}_0)$ is a
``Y'' shaped graph in $\Delta$. Let $\hat{C}= f_0(\hat{C}_0)$. $\{f_0(C_t)\}_{t\in[0,1]}$ is a symplectic isotopy from $f_0(C_0)$ to $\hat{C}$.\\

On the other hand, $\{f_t(C_0)\}_{t\in[0,1]}$ is a symplectic isotopy from $f_0(C_0)$ to $C_0$. Combine the two families we get a symplectic isotopy from $C_0$ to symplectic curve $\hat{C}\subset (\mathbb{C}\hat{\mathbb{P}}^2,\omega_w)$ such that $F_w(\hat{C})=F_0(\hat{C}_0)$ is a ``Y'' shaped graph in $\Delta_w$.
\begin{flushright} $\Box$ \end{flushright}
Now let us get back to the discussion of the mirror of quintic and resume our notation. Based on what we just proved, we can see that there exists a symplectic isotopy from $C = Y_\psi \cap {\rm Sing}(Y_\infty)$ to $\hat{C}$ such that $F_w(\hat{C}) = \Gamma$ is a graph in $\Delta_w$. Use similar argument as in the proof of theorem 9.1 in \cite{lag2}, we have

\begin{theorem}
\label{fe}
There exists a family $\{H_s\}_{s\in [0,1]}$ of piecewise smooth Lipschitz continuous symplectic automorphisms of $\mathbb{C}\hat{\mathbb{P}}^2$ that restrict to identity on the three coordinate $\mathbb{CP}^1$'s, such that $H_0 = {\rm id}$ and $\hat{C} = H_1(C)$.
\end{theorem}
\begin{flushright} $\Box$ \end{flushright}
Let $\Gamma = \Gamma^1\cup \Gamma^2\cup \Gamma^3$, where $\Gamma^1$ is the smooth part of $\Gamma$, $\Gamma^3$ is the singular part of $\Gamma$ in the 1-skeleton of $\Delta_w$, $\Gamma^2$ is the rest of singular part of $\Gamma$. As in the proof of the theorem \ref{be}, we can use theorem \ref{fe} to produce a topologically smooth Lagrangian fibration $\hat{F}_\infty: Y_\infty \rightarrow \partial \Delta_w$ such that $\hat{F}_\infty(C)=\Gamma$. Apply theorem \ref{fd}, we have\\
\begin{theorem}
\label{fb}
For $Y_\psi \subset P_{\Sigma^w}$, when $w=(w_m)_{m\in \Delta^0}\in \tau$ is generic and near the large radius limit of $\tau$, start with Lagrangian fibration $\hat{F}_\infty$, the gradient flow method will produce a Lagrangian fibration $\hat{F}_{\psi}: Y_{\psi} \rightarrow \partial\Delta_w$ with the graph singular locus $\Gamma = \hat{F}_\infty(C)$. There are 4 types of fibers.\\
(i). For $p\in \partial\Delta_w\backslash \Gamma$, $\hat{F}_{\psi}^{-1}(p)$ is a smooth Lagrangian 3-torus.\\
(ii). For $p\in \Gamma^1$, $\hat{F}_{\psi}^{-1}(p)$ is a type I singular fiber.\\
(iii). For $p\in \Gamma^2$, $\hat{F}_{\psi}^{-1}(p)$ is a type II singular fiber.\\
(iv). For $p\in \Gamma^3$, $\hat{F}_{\psi}^{-1}(p)$ is a type III singular fiber.
\end{theorem}
\begin{flushright} $\Box$ \end{flushright}
{\bf Remark:} The advantage of the approach here is that we may start with any \k form $\omega_w$ in the \k class of $w$. The disadvantage is that we have to impose restriction that $F_0(C_0)\subset \Delta_w$ for each curved triangular piece of the singular locus. This condition essentially requires the \k class $w$ to be not too far away from the ray of multiple of the anti-canonical class. To remove such restriction and to be able to construct Lagrangian fibrations with codimension two singular locus for Calabi-Yau hypersurfaces in more general toric varieties, we will use a different method based on a careful construction of the \k form $\omega_w$ so that the argument of deformation to codimension 2 will still work directly under such \k form. More details on this alternative construction of Lagrangian torus fibration of the mirror of quintic Calabi-Yau can be found in \cite{tor}, where the construction of Lagrangian torus fibration and their singular locus and singular fibers are discussed for general Calabi-Yau hypersurfaces in toric varieties.\\

\se{Symplectic topological SYZ construction for quintic Calabi-Yau manifolds}
According to SYZ conjecture, on each Calabi-Yau manifold, there should be a special Lagrangian torus fibration. The mirror Calabi-Yau manifold is the moduli space of special Lagrangian torus together with a flat $U(1)$ bundle on the special Lagrangian torus, with suitable compactification. In particular, the special Lagrangian torus fibration on the mirror should be naturally dual to the special Lagrangian torus fibration on the original Calabi-Yau manifold over the smooth part of the fibration. For our purpose, we will concern the corresponding statement for Lagrangian torus fibrations.\\

After the Lagrangian torus fibrations of a Calabi-Yau manifold and its mirror are constructed, to justify the SYZ conjecture, one first need to find a natural identification of the base spaces (topologically $S^3$) of the Lagrangian fibrations, under which the singular locus of the two Lagrangian fibrations are naturally identified. Then one can compute the monodromy of the two fibration to see if they are dual to each other.\\

\subsection{Identification of the bases}
Apriori the two singular locus arise via very different mechanism and do not seem to match. Nevertheless, there is a natural identification of the base spaces of the two fibrations that naturally identify the singular locus as stated in theorem \ref{fa}. We will carry out the identification of fibers and compute the monodromy in the next subsection.\\

For any $u=(u_m)_{m\in \Delta^0}\in (\tilde{N}\otimes_\mathbb{Z}\mathbb{R} + i\tau)/\tilde{N}$, consider the quintic Calabi-Yau $X_u$ defined by\\
\[
p_u(z)= \sum_{m\in \Delta^0} e^{2\pi iu_m}z^m + z^{m_0}=0.
\]\\
Let $w_m = {\rm Im}(u_m)$, then $w=(w_m)_{m\in \Delta^0}\in \tau = \displaystyle \bigcup_{Z\in \tilde{Z}} \tau_Z$. Assume that $w$ is generic, then $w$ determines a simplicial decomposition $Z$ of $\Delta^0$, and $w\in \tau_Z$. Recall that $Z$ determines a graph $\Gamma_Z \subset \Delta^0$ as the union of simplices in the baricenter subdivision of the simplicial decomposition $Z$ of $\Delta^0$, without any integral points of $\Delta^0$ as vertex. According to theorem \ref{be}, when $w$ is near the large complex limit in $\tau_Z$, (consequently, $X_u$ is in the chamber $U_Z \cong \tau^\mathbb{C}_Z$,) by gradient flow method, we can construct Lagrangian torus fibration of $X_u$ over $\partial \Delta$ such that the singular locus is exactly $\Gamma_Z$.\\

In the mirror side, $w=(w_m)_{m\in \Delta^0}$ can be viewed as representing the \k class of the mirror of quintic. As discussed at the beginning of section 4, $w=(w_m)_{m\in \Delta^0}$ naturally determines a real polyhedron $\Delta_w\subset N$. Theorem \ref{fb} asserts that by gradient flow method, we can construct Lagrangian fibration of $Y_w$ (mirror of quintic) over $\partial \Delta_w$ with the singular locus $\Gamma'_Z$. According to proposition \ref{fc}, The singular locus $\Gamma'_Z$ is naturally the union of 1-simplices in the baricenter subdivision of $\Delta_w$ that under the map $\pi$ do not meet the vertices of $\Delta^\vee$.\\

As a real polyhedron, $\Delta_w$ has a dual real polyhedron $\Delta_w^\vee\subset M$. $\Delta_w^\vee$ has another very simple description, as the convex hull of $\{w_m^{-1}m\}_{m\in \Delta}$. $\{w_m^{-1}m\}_{m\in \Delta\backslash\{m_0\}}$ is exactly the set of vertices of $\Delta_w^\vee$. According to this interpretation, there is a natural piecewise linear map $h: \Delta \rightarrow \Delta_w^\vee$ that maps integral points in $\partial \Delta$ to vertices of $\partial \Delta_w^\vee$. The map is actually a simplicial isomorphism from $\partial \Delta$ to $\partial \Delta_w^\vee$ with respect to the simplicial decomposition $Z$ of $\Delta$. Namely $Z$ is exactly the pullback of the simplicial decomposition of $\Delta_w^\vee$ via the piecewise linear map $\Delta \rightarrow \Delta_w^\vee$ . From this point of view, the Lagrangian fibration of $X_u$ can also be thought of as a Lagrangian fibration over $\Delta_w^\vee$. Let $(\Delta_w^\vee)^0$ denote the image of the 2-skeleton $\Delta^0$ of $\Delta$ into $\Delta_w^\vee$. The image of $\Gamma_Z$, still denoted by $\Gamma_Z$, is the union of simplices in the baricenter subdivision of $(\Delta_w^\vee)^0$ that do not meet vertices of $(\Delta_w^\vee)^0$.\\

Above discussions have reduced the identification of the base spaces of the Lagrangian fibrations for the quintic Calabi-Yau $X_u\in U_Z$ (with the singular locus $\Gamma_Z$) and its mirror $Y_w \in \tau_Z$ (with the singular locus $\Gamma'_Z$) to the problem of identifying $(\Gamma_Z,\partial \Delta_w^\vee)$ and $(\Gamma'_Z, \partial \Delta_w)$.\\

In general, for any convex polyhedron $\Delta$, there is always a natural piecewise linear identification of $\Delta$ and its dual polyhedron $\Delta^\vee$. For any face $\alpha$ of $\Delta$, recall the dual face of $\alpha$ in $\Delta^\vee$ is denoted as $\alpha^*$. Let $\hat{\alpha}$ ($\hat{\alpha}^*$) denote the baricenter of $\alpha$ ($\alpha^*$). Then we have the well known correspondence:

\begin{prop}
The dual correspondence naturally induces a piecewise linear homeomorphism $\partial \Delta \rightarrow \partial \Delta^\vee$ with respect to the baricenter subdivision. Under this homomorphism the baricenter $\hat{\alpha}$ of a face $\alpha$ of $\Delta$ is mapped to the baricenter $\hat{\alpha}^*$ of the dual face $\alpha^*$ on $\Delta^\vee$, and simplex $\{\hat{\alpha}_k\}_{k=0}^l$ is mapped to simplex $\{\hat{\alpha}^*_k\}_{k=0}^l$ linearly.
\end{prop}
\begin{flushright} $\Box$ \end{flushright}
Applying this general result to $\Delta_w$, we get our first topological result on SYZ conjecture.

\begin{theorem}
\label{fa}
There is a natural piecewise linear homeomorphism $\partial \Delta_w \rightarrow \partial \Delta_w^\vee$ that identifies the singular locus $\Gamma'_Z$ and $\Gamma_Z$.\\
\end{theorem}
{\bf Proof:} The two vertices of a 1-simplex in $\Gamma'_Z$ are the baricenters of a 1-dimensional face and a 2-dimensional face of $\Delta_w$ that do not map entirely to a vertex of $\Delta^\vee$ under $\pi$.\\

The two vertices of a 1-simplex in $\Gamma_Z$ are the baricenters of a 1-dimensional simplex and a 2-dimensional simplex in $(\Delta_w^\vee)^0$.\\

Under the dual map, one and two dimensional faces of $\Delta_w$ that do not map entirely to a vertex of $\Delta^\vee$ under $\pi$ are naturally dual to two and one dimensional simplices in $(\Delta_w^\vee)^0$.
\begin{flushright} $\Box$ \end{flushright}

\subsection{Duality of fibers and monodromy}
To fully establish SYZ construction, we also need to establish the duality relation of the Lagrangian torus fibers in the Lagrangian fibrations of quintic Calabi-Yau hypersurfaces and their mirror manifolds. There are many ways to see this, for example, one may compute monodromy operators of the two fibrations, and show that they are dual to each other. We will use a more direct method, we will establish an almost canonical characterization of Lagrangian torus fibers of the Lagrangian fibrations of quintics and the mirror manifolds, and prove that under this almost canonical characterization the two fibers are naturally dual to each other.\\

Consider the Calabi-Yau quintic $X$ defined by the quintic polynomial\\
\[
p_u(z)= \sum_{m\in \Delta^0} e^{2\pi iu_m}z^m + z^{m_0}= \sum_{m\in \Delta^0} e^{2\pi i\eta_m}t^{w_m}z^m + z^{m_0}=0,
\]\\
where $u = \{u_m\}_{m\in \Delta^0}$ and $u_m = \eta_m - i\frac{\log t}{2\pi} w_m$. Lagrangian torus fibration $X\rightarrow \partial \Delta$ is constructed by deforming the natural Lagrangian torus fibration of the large complex limit $X_\infty$ via gradient flow, where $X_\infty$ is defined by\\
\[
p_\infty(z)= z^{m_0}=0.
\]\\
For any vertex $n$ of $\Delta^\vee\subset N$, there is a corresponding 3-dimensional face $\alpha_n$ of $\Delta$ defined as\\
\[
\alpha_n = \{m\in \Delta|\langle m,n\rangle =-1\}.
\]\\
Namely $n$ is the unique supporting function of $\alpha_n$. Clearly, fibers of the fibration $X_\infty\rightarrow \partial \Delta$ over $\alpha_n^0$ (interior of $\alpha_n$) are naturally identified with\\
\[
T_n \cong (N_n \otimes_\mathbb{Z} \mathbb{R})/N_n,
\]\\
where $N_n = N/\{\mathbb{Z}n\}$. Since Lagrangian torus fibration $X\rightarrow \partial \Delta$ is a deformation of fibration $X_\infty\rightarrow \partial \Delta$, we have\\
\begin{prop}
\label{gc}
3-torus fibers of the Lagrangian fibration $X\rightarrow \partial \Delta$ over the interior of $\alpha_n$ can be naturally identified with $T_n$.
\end{prop}
\begin{flushright} $\Box$ \end{flushright}
The dual torus of $T_n$ is naturally\\
\[
T_n^\vee \cong (M_n \otimes_\mathbb{Z} \mathbb{R})/M_n,
\]
where $M_n = n^\perp \subset M$.\\

Similarly, for any integral $m\in \partial\Delta$, we can define\\
\[
T_m^\vee \cong (N_m \otimes_\mathbb{Z} \mathbb{R})/N_m,
\]
where $N_m = m^\perp \subset N$. By multiplication, we have the map $m: T_N^{\mathbb{R}} \rightarrow S^1$. Clearly, $T_m^\vee = m^{-1}(0)$. Similarly, for $\eta_m\in S^1$, we can define $T_{m,\eta_m}^\vee = m^{-1}(\eta_m)$. We have\\
\begin{prop}
\label{ga}
When $X$ is near the large complex limit, 3-torus fiber $X_b$ of the Lagrangian fibration $X\rightarrow \partial \Delta$ over $b$ in a small neighborhood $U_m$ of integral $m\in \partial \Delta$ can be identified with $T_{m,\eta_m}^\vee$. In addition, the identification of $X_b$ and $T_n$ (as discussed in proposition \ref{gc}) can be modified by automorphisms of $T_n$ close to the identity map and smoothly depending on $b\in \alpha_n^0$, such that for $b\in \alpha_n^0 \cap U_m$ the following diagram commutes.\\
\[
\begin{array}{ccc}
X_b&\rightarrow& T_{m,\eta_m}^\vee\\
\downarrow&&\downarrow\\
T_n&\leftarrow& T_N^{\mathbb{R}}\\
\end{array}
\]\\
\end{prop}
{\bf Proof:}
Recall that the moment map\\
\[
F_{w}(x) = \sum_{m\in\Delta} \frac{|x^m|_{w}^2}{|x|_{w}^2}m
\]\\
maps $\mathbb{CP}^4$ to $\Delta$. For integral $m'\in \partial \Delta$, apply propositions 3.2, 3.3, 3.4 in \cite{N} to the simplex $S=\{m_0,m'\}$, we have a $\delta>0$ and a small neighborhood $U_{m'}$ of $m'\in \partial \Delta$ such that

\[
e^{2\pi iu_m}z^m = O(t^{w_{m'} + \delta}),\ \ \ {\rm for}\ m\in \Delta^0\backslash \{m'\}\ {\rm and}\ z\in F_w^{-1}(U_{m'}).
\]

Let $p_u =p_\infty + \tilde{p}_u$, where $p_\infty = z^{m_0}$. When $X$ is near the large complex limit, (namely, $\{w_m\}_{m\in \Delta^0}$ is convex with respect to $\Delta^0$ and $t$ is small,) for $z\in F_w^{-1}(U_{m'})$ near $F_w^{-1}(m')$, we have

\[
\tilde{p}_u(z)= \sum_{m\in \Delta^0} e^{2\pi iu_m}z^m = e^{2\pi iu_{m'}}z^{m'}(1+O(t^{\delta})).
\]

The gradient flow will flow $X = X_1$ to $X_\infty$ through the family $\{X_\psi\}_{\psi \in [1,+\infty)}$, where $X_\psi$ is defined by $p_\psi(z) = \psi z^{m_0} + \tilde{p}_u(z)$. More precisely, we actually use the normalized flow of $V = \frac{\nabla f}{|\nabla f|^2}$, where $f={\rm Re}(q)$ and\\
\[
q= \frac{p_\infty(z)}{\tilde{p}_u(z)} = e^{2\pi iu_{m'}}z^{m_0-m'}\rho(z),\ \ {\rm where}\ \rho(z) = 1+O(t^{\delta}).
\]

The flow produces Lagrangian fibrations $F_\psi: X_\psi \rightarrow \partial \Delta$ from $F_\infty: X_\infty \rightarrow \partial \Delta$.\\

Since $\tilde{p}_u(z)$ is non-vanishing near $F_w^{-1}(m')$, the flow of $V$ moves Lagrangian fibers of $X_\infty$ near $F_w^{-1}(m')$ entirely away from $X_\infty$ to Lagrangian 3-torus fibers of $X_\psi$ inside the large complex torus $T_N$. For monodromy computation, we need to show that these 3-torus fibers in $T_N$ are of the same homotopy class as $T_{m'}^\vee \subset T_N$. Our construction is more precise here, actually identifying these 3-torus fibers with $T_{m,\eta_m}^\vee$.\\

Recall that $T_N = (N \otimes_\mathbb{Z} \mathbb{C})/N$. Let $T_N^{\mathbb{R}} = (N \otimes_\mathbb{Z} \mathbb{R})/N$. Then there is the natural map $P: T_N \rightarrow T_N^{\mathbb{R}}$. For any $b\in U_{m'}$, the fiber $X_b = F_1^{-1}(b)$ deforms through $X_{\psi,b} = F_\psi^{-1}(b)$ according to the inverse flow. Recall from \cite{lag2} that before running the gradient flow, we need to perturb the \k metric to toroidal metric near $C = X_\psi \cap {\rm Sing}(X_\infty)$. In the region we are discussing about, $X_\psi$ will move, therefore is away from $C$ and the \k metric does not need to be modified and is still the original toric metric corresponding to the moment map $F_w$. The local structure of the normalized gradient flow is summarized in theorems 5.3 and 5.4 in \cite{lag2}, which implies that $P(X_{\psi,b})$ will converge to

\[
P(X_{\infty,b}) = \left\{\theta \in T_N^{\mathbb{R}} \left|\langle m'-m_0,\theta \rangle = \eta_{m'} + \frac{{\rm Im}(\log \rho(z))}{2\pi}\right.\right\}
\]

along the inverse flow when $\psi$ approaches $\infty$. In such a way, $X_b$ can be identified with $P(X_{\infty,b}) \subset T_N^{\mathbb{R}}$. (The identification in proposition \ref{gc} also factor through this identification.)\\

Notice that

\[
T_{m',\eta_{m'}}^\vee = \left\{\theta \in T_N^{\mathbb{R}} \left|\langle m'-m_0,\theta \rangle = \eta_{m'}\right.\right\}
\]

and $\log \rho(z) = O(t^{\delta})$ can be arbitrarily small when $X$ is near the large complex limit. Namely, $P(X_{\infty,b})$ are small perturbations of $T_{m',\eta_{m'}}^\vee$ that smoothly depend on $b\in U_{m'}$. A projection along some direction transverse to $T_{m',\eta_{m'}}^\vee$ will identify all these $P(X_{\infty,b})$ with $T_{m',\eta_{m'}}^\vee$ smoothly depending on $b\in U_{m'}$. In this way, we get the desired identifications $X_b \rightarrow T_{m',\eta_{m'}}^\vee$ for $b\in U_{m'}$.\\

Since both identifications $X_b \rightarrow T_{m',\eta_{m'}}^\vee$ and $X_b \rightarrow T_n$ factor through $P(X_{\infty,b})$, the diagram in the proposition is equivalent to the following diagram\\
\[
\begin{array}{ccc}
P(X_{\infty,b})&\rightarrow& T_{m,\eta_m}^\vee\\
\downarrow&&\downarrow\\
T_n&\leftarrow& T_N^{\mathbb{R}}\\
\end{array}
\]\\

Since $P(X_{\infty,b})$ is a small perturbation of $T_{m,\eta_m}^\vee$, and the map $P(X_{\infty,b}) \rightarrow T_n$ also factor through $T_N^{\mathbb{R}}$ naturally, the automorphism of $T_n$ defined as composition of the 4 arrows in the diagram (reversing the arrow of the map $P(X_{\infty,b}) \rightarrow T_n$) will be as desired and is close to the identity map of $T_n$.
\begin{flushright} $\Box$ \end{flushright}
Now, let us consider the mirror Calabi-Yau $Y\subset P_{\Sigma^w}$. Also by gradient flow method, we can construct Lagrangian fibration $Y \rightarrow \partial \Delta_w$. For any $m \in \partial \Delta \cap M$, there is a corresponding 3-dimensional face $\alpha_m$ of $\Delta_w$ defined as\\
\[
\alpha_m = \{n\in \Delta_w|w_m^{-1}\langle m,n\rangle =-1\}.
\]\\
Namely $w_m^{-1}m$ is the unique supporting function of $\alpha_m$. Clearly, fibers of the fibration $Y_\infty\rightarrow \partial \Delta_w$ over $\alpha_m^0$ (interior of $\alpha_m$) are naturally identified with\\
\[
T_m \cong (M_m \otimes_\mathbb{Z} \mathbb{R})/M_m,
\]\\
where $M_m = M/\{\mathbb{Z}m\}$. Since Lagrangian torus fibration $Y\rightarrow \partial \Delta_w$ is a deformation of fibration $Y_\infty\rightarrow \partial \Delta_w$, we have\\
\begin{prop}
3-torus fibers of the Lagrangian fibration $Y \rightarrow \partial \Delta_w$ over interior of $\alpha_m$ can be naturally identified with $T_m$.
\end{prop}
\begin{flushright} $\Box$ \end{flushright}
Recall the natural map $\pi: \partial \Delta_w \rightarrow \partial \Delta^\vee$. For any vertex of $n$ of $\Delta^\vee$, we have\\
\begin{prop}
\label{gb}
3-torus fiber $Y_b$ of the Lagrangian fibration $Y \rightarrow \partial \Delta_w$ over $b$ in a small neighborhood $U_n$ of $\pi^{-1}(n) \subset \partial \Delta_w$ can be naturally identified with $T_n^\vee$. In addition, if $b\in \alpha_m^0$, then the following diagram commutes.\\
\[
\begin{array}{ccc}
Y_b&\rightarrow& T_n^\vee\\
\downarrow&&\downarrow\\
T_m&\leftarrow& T_M^\mathbb{R}
\end{array}
\]
\end{prop}
\begin{flushright} $\Box$ \end{flushright}
This proposition can be proved in a way similar to proposition \ref{ga}. For detail, please see \cite{tor}, where a more general result is proved.\\

Recall from the theorem \ref{fa}, we have the natural piecewise linear identification of the two base of the Lagrangian fibrations\\
\[
s: \partial \Delta_w \rightarrow \partial \Delta_w^\vee \cong \partial \Delta.
\]\\
Let $\alpha_n^0$, $\alpha_m^0$ denote the interior of $\alpha_n$, $\alpha_m$. Then we have\\
\begin{prop}
$U_m$, $U_n$ can be suitably chosen such that\\
\[
s(U_n) = \alpha_n^0,\ \ s(\alpha_m^0) = U_m.
\]
And\\
\[
\partial \Delta_w \backslash \Gamma' = \left(\bigcup_{m\in \Delta^0} \alpha_m^0\right) \cup \left(\bigcup_{n\in \partial \Delta^\vee} U_n\right),
\]\\
\[
\partial \Delta \backslash \Gamma = \left(\bigcup_{n\in \partial \Delta^\vee} \alpha_n^0\right) \cup \left(\bigcup_{m\in \Delta^0} U_m\right).
\]
\end{prop}
{\bf Proof:} When we enlarge $U_m$, we need to extend the identification of the fibers with $T_{m,\eta_m}^\vee$ and modify identification of fibers over $\alpha_n^0$ with $T_n$ accordingly to make sure that the diagram still commute (as in proposition \ref{ga}). Such operations are very easy to do, because $U_m$ is contractible and the fibration of $X_\psi$ over $U_m$ is trivial and smooth.
\begin{flushright} $\Box$ \end{flushright}
We are now ready to establish the dual relation of the fibers. Our construction of the torus fibration for $Y$ only depends on the real \k moduli, which under the monomial-divisor map corresponds to the moduli of quintics with real coefficients. Therefore, we will restrict our duality discussion to quintic $X$ in such real moduli space, where duality relation is much more precise and without shift. Notice that when the quintic has real coefficients, we have $\eta_m=0$ for all $m$ and $T_{m,\eta_m}^\vee = T_m^\vee$. For any $b\in \partial \Delta_w \backslash \Gamma'$, let $Y_b$ ($X_{s(b)}$) denote the fiber of the Lagrangian fibration $Y \rightarrow \partial \Delta_w$ ($X \rightarrow \partial \Delta$) over $b\in \partial \Delta_w \backslash \Gamma'$ ($s(b)\in \partial \Delta \backslash \Gamma$). Then we have\\
\begin{theorem}
For any $b\in \partial \Delta_w \backslash \Gamma'$, $Y_b$ is naturally dual to $X_{s(b)}$.\\
\end{theorem}
{\bf Proof:}
Based on above propositions, duality is easy to establish. The only thing that need to be addressed is that duality defined in two ways according to $U_m$ or $\alpha_n^0$ for $b\in U_m\cap \alpha_n^0$ coincide. For this purpose, one only need to show that the following diagrams commute.\\
\[
\begin{array}{ccccccc}
X_{s(b)}&\rightarrow& T_m^\vee & \ \ &Y_b&\rightarrow& T_n^\vee \\
\downarrow&&\downarrow& \ &\downarrow&&\downarrow\\
T_n&\leftarrow& T_N^\mathbb{R}& \ &T_m&\leftarrow& T_M^\mathbb{R}\\
\end{array}
\]\\
This is proved in proposition \ref{ga} and \ref{gb}.
\begin{flushright} $\Box$ \end{flushright}
With the explicit identification of fibers in place, monodromy computation becomes a piece of cake! Consider the path $\gamma_{nmn'm'} = \alpha_n^0 U_m \alpha_{n'}^0 U_{m'} \alpha_n^0$ on $\partial \Delta$, where $n,n'\in \partial \Delta^\vee$, $m,m'\in \Delta^0$ satisfy\\
\[
\langle m,n\rangle=\langle m,n'\rangle = \langle m',n\rangle = \langle m',n'\rangle = -1.
\]\\
This condition implies that $\alpha_n$ and $\alpha_{n'}$ have common face that contains $m,m'$. Correspondingly we have the diagram\\
\[
\begin{array}{ccc}
N_n&\rightarrow& N_m\\
\uparrow&&\downarrow\\
N_{m'}&\leftarrow& N_{n'}\\
\end{array}
\]\\
\[
x\rightarrow x+\langle m,x\rangle n \rightarrow x+\langle m,x\rangle n + \langle m',x+\langle m,x\rangle n\rangle n'= x+\langle m,x\rangle n + \langle m'-m,x\rangle n'
\]\\
Compose the four operators and modulo $n$, we get\\
\begin{theorem}
The monodromy operator along $\gamma_{nmn'm'}$ is\\
\[
[x] \rightarrow [x] + \langle m'-m,x\rangle [n']\ \ \ {\rm for}\ [x]\in N_n.
\]
\end{theorem}
\begin{flushright} $\Box$ \end{flushright}
{\bf Remark:} Now we have found an extremely simple way to compute monodromy. All our monodromy computation in \cite{lag1} can be much more easily performed by this method. Since monodromy computation is becoming so trivial, I will omit the corresponding computation of monodromy operator for the mirror fibration, which is naturally dual to the monodromy operator for the fibration of quintic.\\\\

\se{Singular fibers of Lagrangian torus fibration}
In the previous section we established the duality between the smooth torus fibers of the Lagrangian torus fibrations for quintics and their mirrors. In this section we will discuss the generic singular fibers and their duality.\\

Let us start with some facts for $C^\infty$ Lagrangian fibrations. Let $(X,\omega)$ be a smooth symplectic manifold. A fibration $F: X\rightarrow B$ is called a $C^l$-Lagrangian fibration, if $F$ is a $C^l$ map and the smooth part of each fiber is \l and belongs to the regular point set of map $F$. The following well known results were discussed in the section 2 of \cite{lag2}.\\
\begin{theorem}
Let $F: X\rightarrow B$ be a $C^{1,1}$-Lagrangian fibration, then for any $b\in B$, there is an action of $T^*_bB$ on $X_b = f^{-1}(b)$.
\end{theorem}
\begin{flushright} $\Box$ \end{flushright}
{\bf Remark:} In particular, this theorem applies to $C^\infty$ Lagrangian fibrations.\\

\begin{co}
\label{ea}
For any $b \in B$,\\
\[
{\rm Reg}(X_b) = {\rm Reg}(F^{-1}(b)) = \bigcup_l O_l
\]
is a disjoint union of orbits of $T^*_bB$, where each $O_l$ is diffeomorphic to $(S^1)^k\times\mathbb{R}^m$ for some $k+m = \dim B$.
\end{co}
\begin{flushright} $\Box$ \end{flushright}
{\bf Remark:} When $F$ is generic in certain sense, we expect $\displaystyle F^{-1}(b) = \bigcup_l O_l$ to be a disjoint union of finitely many orbits of $T^*_bB$, where each $O_l$ is diffeomorphic to $(S^1)^k\times\mathbb{R}^m$ for some $k+m \leq \dim B$.\\

From these facts, one can see that $C^\infty$ Lagrangian fibrations put rather strict restriction on the topology of singular fibers. In the following, we will discuss certain generic singular fibers that include the singular fibers which appeared in the Lagrangian torus fibrations we constructed for generic quintics and Fermat type quintics. Although these singular fibers all satisfy the topological constraints in corollary \ref{ea}, it is not immediately clear whether they can be realized as singular fibers for $C^\infty$ Lagrangian torus fibrations, especially the type two singular fibers.\\

As we know, even for elliptic fibrations the singular fibers in general can be quite complicated. It is crucial to restrict our attention to some classes of stable singular fibers with certain generic nature.\\

In three dimensions, singular fibers conceivably can be even more complicated. To have a meaningful discussion, it is crucial for us to first concentrate on certain classes of ``generic" singular fibers. We will restrict our discussion to three types of singular fibers.\\

Type one singular fiber comes from the product of a 2-dimensional singular fiber with a circle. It has one vanishing 1-cycle. In particular, we denote the product of 2-dimensional $A_n$ singular fiber with a circle by $I_n$.\\

Type two singular fiber has one vanishing 1-cycle, and has a natural map to a 2-torus with fiber being either a point or a circle representing the vanishing cycle.\\

Type three singular fiber has two independent vanishing 1-cycles, and has a natural map to a circle with fiber being either a point or a 2-torus representing the vanishing cycles.\\

\begin{prop}
Type three singular fibers are parametrized by positive integers, denoted by $III_n$. $n$ is the number of point fibers for the map to the circle. $III_n$ fiber has Euler number equal to $n$. In particular, type $III(=III_1)$ fiber is the generic singular fiber with Euler number 1.
\end{prop}
\begin{flushright} $\Box$ \end{flushright}
\begin{center}
\leavevmode
\hbox{%
\epsfxsize=4in
\epsffile{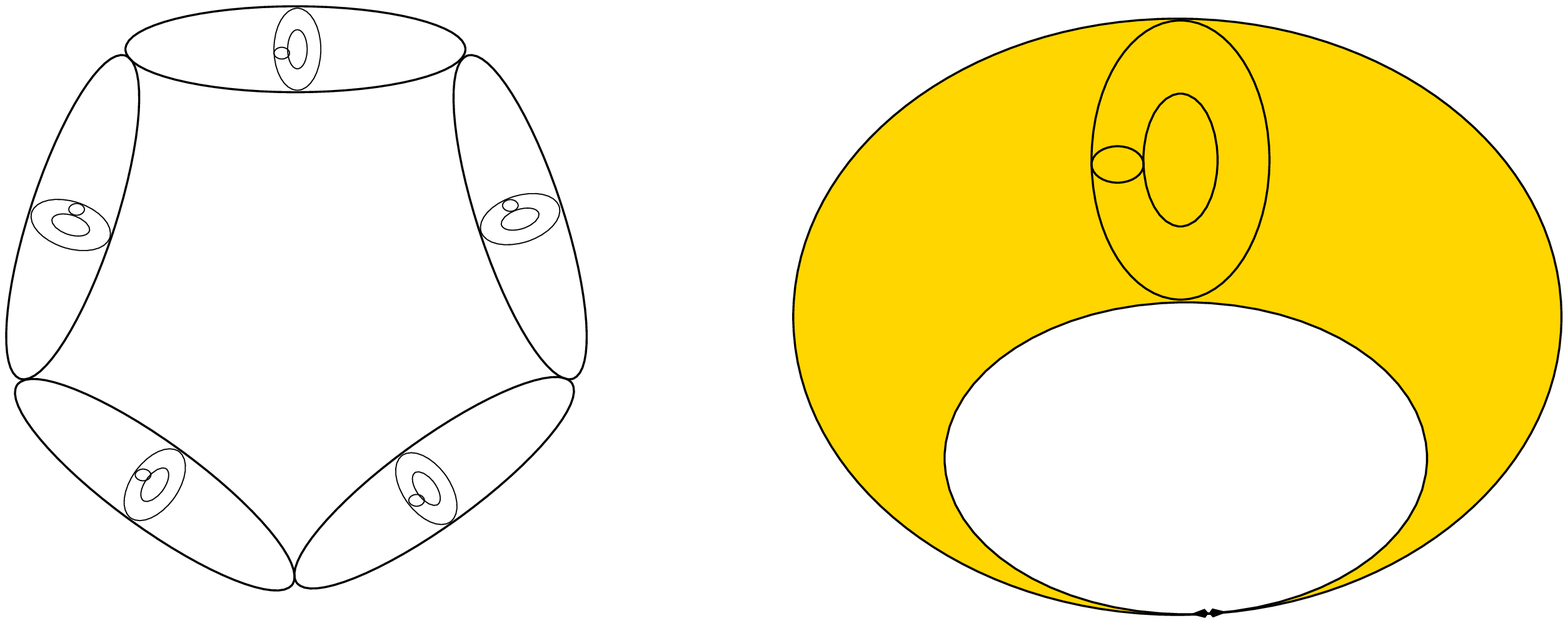}}
\end{center}
\begin{center}
\stepcounter{figure}
Figure \thefigure: Type $III_5$ and type $III$ fibers
\vspace{10pt}
\end{center}
Type two singular fiber has a map to a 2-torus $T^2$. The set of points on the 2-torus with point fiber is typically a graph $\Gamma$ on the 2-torus, which divides the 2-torus into several regions. In general a graph $\Gamma$ that divides $T^2$ into $n$ region has Euler number equal to $n$. We have\\
\begin{prop}
For a type two singular fiber, if the corresponding graph $\Gamma \subset T^2$ devides the 2-torus $T^2$ into $n$ regions, then the Euler number for the singular fiber is equal to $-n$.
\end{prop}
\begin{flushright} $\Box$ \end{flushright}
We are interested in the generic type two singular fiber with Euler number equal to $-1$. We have\\
\begin{prop}
There are only two type two singular fibers with Euler number equal to $-1$, the type $II$ fiber corresponding to parallel sexgon and type $\tilde{II}$ fiber corresponding to parallelgram.\\
\end{prop}
{\bf Proof:}
Since Euler number equals to $-1$, by the previous proposition, $\Gamma$ divides $T^2$ into only one region. Going to the universal cover of $T^2$, this one region gives a filling of $\mathbb{R}^2$ by only one type of polygon. It is well known that there are only two types of polygon that can fill the plane, the parallel sexgon and the parallelgram.
\begin{flushright} $\Box$ \end{flushright}
\begin{center}
\leavevmode
\hbox{%
\epsfxsize=4in
\epsffile{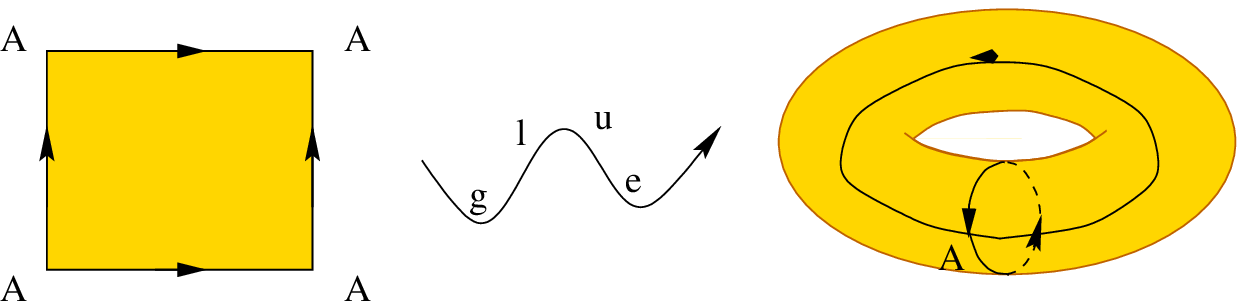}}
\end{center}
\begin{center}
\stepcounter{figure}
Figure \thefigure: Type $\tilde{II}$ fiber
\end{center}
\vspace{10pt}
\begin{center}
\leavevmode
\hbox{%
\epsfxsize=4in
\epsffile{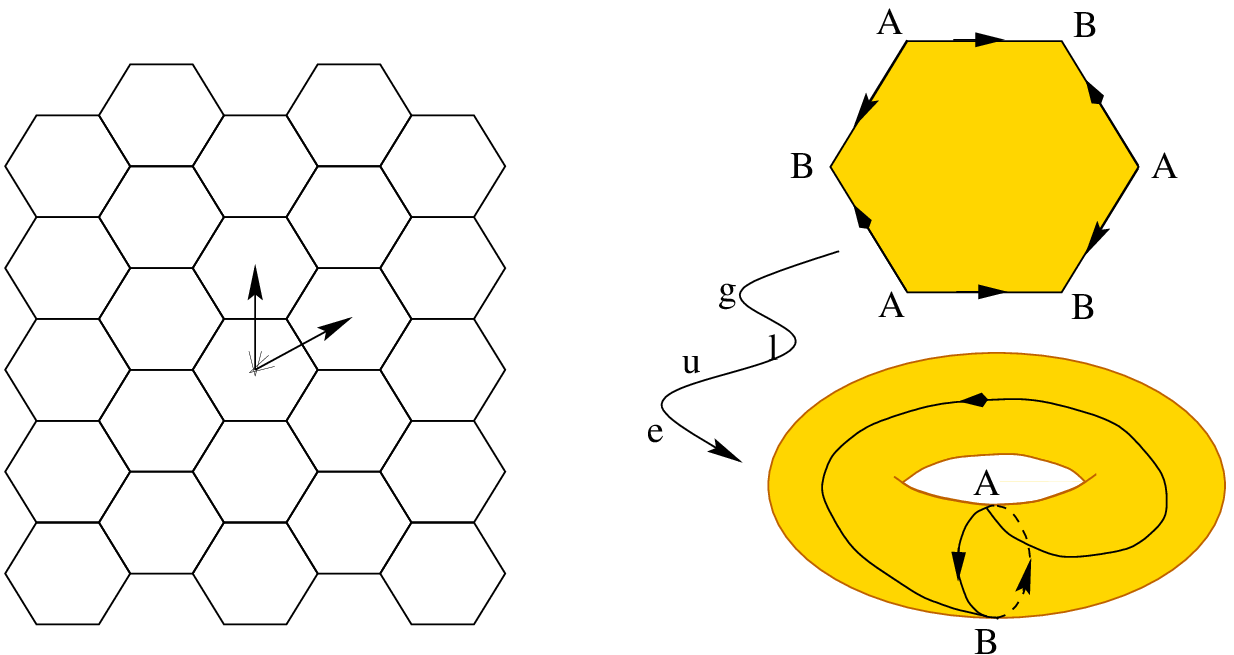}}
\end{center}
\begin{center}
\stepcounter{figure}
Figure \thefigure: Type $II_{5\times 5}$ and type $II$ fibers
\vspace{10pt}
\end{center}
$n\times m$ cover of type $II$ ($\tilde{II}$) singular fiber are type two singular fiber with Euler number equal to $mn$. We will denote this kind of singlar fiber by $II_{n\times m}$ ($\tilde{II}_{n\times m}$).\\

In our generic Lagrangian torus fibrations for quintics, we have generic singular fiber type $I$, $II$, $III$. Type $I$ fiber is dual to type $I$ fiber, type $II$ fiber is dual to type $III$ fiber. Clearly, the dual relations change only the sign of the Euler number of the singular fiber, while keeping the absolute value unchanged. This apparently is in accordance with mirror symmetry.\\

In our discussion of Lagrangian torus fibrations of Fermat type quintics, we have singular fibers of type $I_5$, $II_{5\times 5}$, $III_5$.\\

Summarize our results, we have proved the symplectic topological version of SYZ conjecture for quintic Calabi-Yau hypersurfaces, including the refinement involving the singular fibers. Recall that $\tilde{Z}$ denotes the set of integral simplicial decompositions of $\Delta^0$ (the 2-skeleton of $\Delta$).\\
\begin{theorem}
\label{eb}
For $Z \in \tilde{Z}$, consider generic quintic Calabi-Yau hypersurface $X$ in the chamber $U_Z$ near the large complex limit and its mirror Calabi-Yau manifold $Y$ with the \k moduli $w\in \tau_Z$ near the large radius limit, there exist corresponding Lagrangian torus fibrations\\
\[
\begin{array}{ccccccc}
X_{s(b)}&\hookrightarrow& X& \ \ &Y_b&\hookrightarrow& Y\\
&&\downarrow& \ &&&\downarrow\\
&& \partial \Delta& \ &&& \partial \Delta_w\\
\end{array}
\]\\
with singular locus $\Gamma_Z \subset \partial \Delta$ and $\Gamma'_Z \subset \partial \Delta_w$, where $s:\partial \Delta_w \rightarrow \partial \Delta$ is a natural homeomorphism and $s(\Gamma'_Z)=\Gamma_Z$. For $b\in \partial \Delta_w \backslash \Gamma'_Z$, the corresponding fibers $X_{s(b)}$ and $Y_b$ are naturally dual to each other.\\

Singular fibers over the smooth part of $\Gamma_Z$ ($\Gamma'_Z$) are of type I. Singular fibers over the vertices of $\Gamma_Z$ ($\Gamma'_Z$) are of type II or III. Type I fiber is dual to type $I$ fiber, type II fiber and type III fiber are dual to each other.
\end{theorem}
\begin{flushright} $\Box$ \end{flushright}
The original SYZ mirror conjecture was rather sketchy in nature, with no mentioning of singular locus, singular fibers and duality of singular fibers, which is essential if one wants, for example, to use SYZ to construct mirror manifold. Our construction of Lagrangian torus fibrations and proof of symplectic topological SYZ for generic Calabi-Yau quintic hypersurfaces explicitly produce the 3 types of generic singular fibers (type $I$, $II$, $III$) and determine the way they are dual to each other under mirror symmetry. Our construction clearly indicates what should happen in general. In particular, it suggests that type II singular fiber with Euler number $-1$ should be dual to type III singular fiber with Euler number $1$. This together with the knowledge of singular locus from our construction enable us to give a more precise formulation of SYZ mirror conjecture (in symplectic category). This precise formulation naturally suggests a way to construct mirror manifold in general from a generic Lagrangian torus fibration of a Calabi-Yau manifold.\\

{\bf Precise symplectic SYZ mirror conjecture}
For any Calabi-Yau 3-fold $X$, with Calabi-Yau metric $\omega_g$ and holomorphic volume form $\Omega$, there exists a Lagrangian torus fibration of $X$ over $S^3$\\
\[
\begin{array}{ccc}
T^3&\hookrightarrow& X\\
&&\downarrow\\
&& S^3
\end{array}
\]\\
with a Lagrangian section and codimension 2 singular locus $\Gamma\subset S^3$, such that general fibers (over $S^3\backslash \Gamma$) are 3-torus. For generic such fibration, $\Gamma$ is a graph with only 3-valent vertices. Let $\Gamma = \Gamma^1\cup \Gamma^2 \cup \Gamma^3$, where $\Gamma^1$ is the smooth part of $\Gamma$, $\Gamma^2 \cup \Gamma^3$ is the set of the vertices of $\Gamma$. For any leg $\gamma\subset \Gamma^1$, the monodromy of $H_1(X_b)$ of fiber under suitable basis is\\
\[
T_\gamma =\left(
\begin{array}{ccc}
1 &  1 & 0\\
0 &  1 & 0\\
0 &  0 & 1
\end{array}
\right).
\]\\
Singular fiber along $\gamma$ is of type $I$.\\

Consider a vertex $P\in\Gamma^2 \cup \Gamma^3$ with legs $\gamma_1$, $\gamma_2$, $\gamma_3$. Correspondingly, we have monodromy operators $T_1$, $T_2$, $T_3$.\\

For $P\in\Gamma^2$, under suitable basis we have\\
\[
T_1 =\left(
\begin{array}{ccc}
1 &  1 & 0\\
0 &  1 & 0\\
0 &  0 & 1
\end{array}
\right),\ \
T_2 =\left(
\begin{array}{ccc}
1 &  0 & -1\\
0 &  1 & 0\\
0 &  0 & 1
\end{array}
\right),\ \
T_3 =\left(
\begin{array}{ccc}
1 &  -1 & 1\\
0 &  1 & 0\\
0 &  0 & 1
\end{array}
\right).
\]\\
Singular fiber over $P$ is of type $II$.\\

For $P\in\Gamma^3$, under suitable basis we have\\
\[
T_1 =\left(
\begin{array}{ccc}
1 &  0 & 0\\
1 &  1 & 0\\
0 &  0 & 1
\end{array}
\right),\ \
T_2 =\left(
\begin{array}{ccc}
1 &  0 & 0\\
0 &  1 & 0\\
-1 &  0 & 1
\end{array}
\right),\ \
T_3 =\left(
\begin{array}{ccc}
1 &  0 & 0\\
-1 &  1 & 0\\
1 &  0 & 1
\end{array}
\right).
\]\\
Singular fiber over $P$ is of type $III$.\\

The Lagrangian fibration for the mirror Calabi-Yau manifold $Y$ has the same base $S^3$ and singular locus $\Gamma \subset S^3$. For $b\in S^3\backslash \Gamma$, $Y_b$ is the dual torus of $X_b\cong T^3$. In another word, the $T^3$-fibrations\\
\[
\begin{array}{ccc}
T^3&\hookrightarrow& X\\
&&\downarrow\\
&& S^3\backslash\Gamma
\end{array}
\ \ \ \ \ \
\begin{array}{ccc}
T^3&\hookrightarrow& Y\\
&&\downarrow\\
&& S^3\backslash\Gamma
\end{array}
\]\\
are dual to each other. In particular the monodromy operators will be dual to each other.\\

For the fibration of $Y$, singular fibers over $\Gamma^1$ should be type $I$, singular fibers over $\Gamma^2$ should be type $III$, singular fibers over $\Gamma^2$ should be type $II$. Namely, dual singular fiber of a type $I$ singular fiber is still type $I$. Type $II$ and $III$ singular fibers are dual to each other.\\

{\bf Conjecture:} Type $II$, $\tilde{II}$, $III$ singular fibers are the only possible generic singular fibers with Euler number equal to $\pm 1$.\\

{\bf Remark:} The last piece of the SYZ puzzle we have not yet discussed in detail is the construction of a section of the Lagrangian fibration. With the explicit description of Lagrangian fibers in section 5, it is not hard to construct the section. When the coefficients of the quintic are all real, one can take the identity section on each piece with explicit description as in proposition \ref{ga}. they piece together to form an global section that can be extended to be over the compliment of the singular locus. With a little more care, such section can be made Lagrangian. We will describe the precise construction of global sections in \cite{mon}, where more general Lagrangian cycles in Calabi-Yau hypersurfaces will also be disscussed.\\\\

{\bf Acknowledgement:} I would like to thank Qin Jing for many very stimulating discussions during the course of my work, and helpful suggestions while carefully reading my early draft. I would also like to thank Prof. S.-T. Yau for his constant encouragement. This work was originally done while I was in Columbia University. I am very grateful to Columbia University for excellent research environment.\\\\

\ifx\undefined\bysame
\newcommand{\bysame}{\leavevmode\hbox to3em{\hrulefill}\,}
\fi

\end{document}